\documentclass[a4paper, 12pt]{article}
\usepackage[utf8]{inputenc}
\usepackage{makeidx}
\usepackage{enumerate, theorem}
\usepackage{amsmath, amsfonts, amssymb}
\usepackage[height=22.5cm, width=15cm]{geometry}
\usepackage[all]{xy}
\usepackage{mathrsfs, yfonts}
\usepackage{graphicx}
\usepackage{cancel}
\DeclareMathOperator{\C}{\mathbb{C}}
\newcommand{\parag}[1]{\paragraph{\sc{#1.}}}

\newtheorem{thm}{Theorem}[subsection]
\newtheorem{defn}[thm]{Definition}
\newtheorem{cor}[thm]{Corollary}
\newtheorem{prop}[thm]{Proposition}
\newtheorem{lemma}[thm]{Lemma}

\begin{document}

\title{The $\mu$-Trace System}

 \author{Daniel Barlet\footnote{Institut Elie Cartan, G\'eom\`{e}trie,\newline
Universit\'e de Lorraine, CNRS UMR 7502   and  Institut Universitaire de France.}.}

\maketitle

\hfill  {\it Qu'avez-vous donc trouvé qui vous rend si content? } \newline \smallskip
 \hfill {\it Je ne sais pas très bien,  mais c'est  intéressant !}
 
 \parag{Abstract} We study a simple  1-parameter perturbation of the regular holonomic  Trace System satisfied by a complex  power of the root of the universal polynomial of degree k as a holomorphic function of the coefficients. We prove that these systems have many analogous properties than the Trace System studied in \cite{[B.24]} and we prove that they are, in general, minimal extensions of a simple pole meromorphic connection on a  rank $k$ trivial bundle on $\mathbb{C}^k$. We also examine the structure of these $D$-modules for the special values of the parameters. This explicites many examples of perverse sheaves associated to representations of the $\pi_1$ of the complement of  the hyper-surface $\{\sigma_k\Delta(\sigma) = 0\}$ in the affine space with coordinates $\sigma_1, \dots, \sigma_k$, where $\Delta(\sigma)$ is the discriminant of the universal monic  polynomial  of degree $k$, $P_\sigma(z) := z^k + \sum_{h=1}^k (-1)^h\sigma_h z^{k-h}$.
 
 \parag{Classification} 32 C 38 - 35 Q 15 - 14 F 10 - 35 C 10 - 32 S 40.
 
  \parag{Key Words} Regular holonomic D-module. Minimal extension. Trace System. Perverse sheaf.

\tableofcontents

.

\section{Introduction}

After a detailed study of the "Trace System" (see \cite{[B-MF]}, \cite{[B.21]}  \cite{[B.22]} and \cite{[B.24]})  I was intrigued by the fact that a partial derivative of a  trace function
is solution of a partial differential system obtained by a very simple modification of the Trace System. Precisely the Trace System corresponds to the left ideal $\mathcal{I}$ in the Weyl algebra $\mathbb{C} \langle \sigma_1, \cdots, \sigma_k, \partial_1, \cdots, \partial_k \rangle $ of $N := \mathbb{C}^k$ generated by the $(2, 2)$-minors of the matrix
\begin{equation*}
M := 
\begin{pmatrix}
 \partial_1 & -E \\ 
 \partial_2 & \partial_1 \\ 
 \cdot & \partial_2 \\
  \cdot & \cdot \\ 
  \cdot & \cdot \\
   \partial_k & \partial_{k-1} 
 \end{pmatrix} \quad  where \  E := \sum_{h=1}^k \sigma_h\partial_h.
  \end{equation*}
 So explicitly by the  $ A_{p,q} = \partial_p\partial_q - \partial_{p+1}\partial_{q-1} \quad \forall p \in [1, k-1],  \forall q \in [2, k] $   and also  by the  
    $ T^m := \partial_1\partial_{m-1} + \partial_m E \quad  \forall m \in [2, k]$.  Then the partial derivatives of order $p \geq 1$ of any  trace function (which is, by definition, a solution of the system above) are solutions of the system obtained by changing only  $E$ to $E + p$.\\
   So it is tempting to study the system obtained by changing $E$ to $E + \mu$ with $\mu \in \mathbb{C} \setminus \mathbb{Z}$. I call this system the "$\mu$-Trace System" and I denote by $\mathcal{I}_\mu$ the corresponding ideal of $D_N$, and by $\mathcal{M}_\mu := D_N\big/\mathcal{I}_\mu$ the corresponding $D_N$-module.\\
    
    Several remarkable properties of the Trace System have their analogs for the $\mu$-Trace System:
    \begin{enumerate}
    \item The change from $0$ to $\mu$  does not affect the characteristic variety $Z$ of the system since it is a reduced  and irreducible  $(k+1)$-dimensional variety in $T^*_N$ 
    which has only  $1$-dimensional fibers over $N := \C^k$ (see \cite{[B.22]}).
    \item There exists an action of the Lie algebra $sl_2(\mathbb{Z})$ on the $D_N$-module $\mathcal{M}_\mu$ given by right product by the following  order $1$ differential operators
    \begin{itemize}
    \item $U_{-1} := \sum_{h=0}^{k-1} (k-h)\sigma_h\partial_{h+1} $ with the convention $\sigma_0 \equiv 1$.
    \item $2U_0 + k\mu := 2\sum_{h=1}^k h\sigma_h\partial_h  + k\mu $ 
    \item $U_1 + \mu\sigma_1 :=  \big(\sum_{h=1}^k (\sigma_1\sigma_h - (h+1)\sigma_{h+1})\partial_h\big)  + \mu\sigma_1 $ \ where $\sigma_{k+1} \equiv 0$
    \end{itemize}  satisfying the commutation relations
    \begin{align*}
    & [2U_0+k\mu, U_{-1}] = -2U_{-1},\\
    &  [2U_0 +k\mu, U_1 + \mu\sigma_1] = 2(U_1 + \mu\sigma_1) \quad {\rm and}\\
    &  [U_{-1}, U_1+ \mu\sigma_1] = 2U_0 + k\mu . \tag{@}
    \end{align*}
    \item The left ideal of the $\mu$-Trace System is the annihilator of a sequence of polynomials in $\mathbb{C}[\sigma]$:\\
    -- For $\mu = 0$ the Newton  polynomials $N_m, m \in \mathbb{N}$.\\
    -- For $\mu \in \mathbb{C}\setminus \mathbb{Z}$ the polynomials $\Phi_{\mu, m}(\sigma)$ given, for $R$ large enough compared to $\vert\vert \sigma\vert\vert$,  by the formula
    \begin{equation}
    \Phi_{\mu, m}(\sigma) := \frac{1}{2i\pi}\int_{\vert a\vert = R} a^m\Big(\frac{a^k}{P_\sigma(a)}\Big)^\mu \ \frac{da}{a} ,
    \end{equation}
    where we denote 
    $$P_\sigma(a) := a^k + \sum_{h=1}^k (-1)^h \sigma_ha^{k-h}.$$
    This formula is the analog of the formula:
    \begin{equation}
    N_m(\sigma) = \frac{m}{2i\pi}\int_{\vert a\vert = R}  a^m Log\big(\frac{a^k}{P_\sigma(a)}\big) \frac{da}{a}\quad\quad  \forall m \in \mathbb{N}^*
    \end{equation}
    -- The analog of the formula
    \begin{equation}
   N_m =  (m-1)! \, U_1^{m-1}(N_1)
    \end{equation}
    is given by :
    \begin{equation}
    \Phi_{\mu, m} = m!\, (U_1 + \mu\sigma_1)^{m-1}( \Phi_{\mu, 1}) \quad {\rm where} \quad  \Phi_{\mu, 1} = \mu N_1 = \mu\sigma_1.
    \end{equation}
    For $\mu = 0$ the ideal $\mathcal{I}$ of the Trace System is the annihilator in $D_N$ of the function $Log\big(\frac{a^k}{P_\sigma(a)}\big)$ which is holomorphic
    on an open subset of $N \times \mathbb{C}$ on which $\vert a\vert$ is large enough compared with $\vert\vert\sigma\vert\vert$. For $\mu \in \mathbb{C} \setminus \mathbb{Z}$ the ideal $\mathcal{I}_\mu$ of the $\mu$-Trace System is the annihilator in $D_N$ of the function $\big(a^k/P_\sigma(a))^\mu$ which is holomorphic
    on an open subset of $N \times \mathbb{C}$ on which $\vert a\vert$ is large enough compared with $\vert\vert\sigma\vert\vert$.
    \item As in the case of the Trace System, for $\mu \not\in \mathbb{Z}$ the right multiplication by $\partial_h, h \in [1, k]$ sends the left ideal $\mathcal{I}_{\mu+1}$ to 
    $\mathcal{I}_{\mu}$.\\
    This is also true for the right multiplication by $E + \mu$.  This corresponds to the formulas
    \begin{align}
   & \partial_h \Phi_{\mu, m} = (-1)^{h-1}\mu\Phi_{\mu+1, m-h} \qquad \forall h \in [1, k] \\
   & (E + \mu)[\Phi_{\mu, m}] = \mu\Phi_{\mu+1, m}
    \end{align}
    analogous to the formula for $\mu = 0$:
    \begin{align}
    & \partial_h N_m/m = (-1)^h  M_{m-h} \quad {\rm where} \quad  M_p(\sigma) := \sum_{j=1}^k \frac{z_j^{p+k-1}}{P'(\sigma, z_j)} = \Phi_{1, p}(\sigma)\\
    & E[N_m/m] = M_m = \Phi_{1, m}.
    \end{align}

    \item When the fonction $Log(\sigma_k)$ is well defined on an open subset $U$ of $N$ it is a trace function on $U$.\\
    More generally, for any given $a \in \mathbb{C}^*$, when the fonction $Log\big(a^k/P_\sigma(a)\big)$ is well defined and holomorphic  on an open subset $D$  in  $N $,  it is a trace function.\\
    In fact the  function $Log \big(\frac{a^k}{P_\sigma(a)}\big)$ generates in some sense the trace functions and the trace system because the coefficient of $a^{-m}$ in its Laurent expansion at infinity (in $a$)  is $N_m/m$ and the left ideal $\mathcal{I}$ may be defined as the annihilator of this function (and the constants), or equivalently, as the annihilator  of all the $N_m, m \in \mathbb{N}$.\\
   This corresponds for $\mu \in \C \setminus \mathbb{Z}$  to the following   facts: \\
   - When the function $\sigma_k^{-\mu}$ is defined on an open subset in $N$ it is a $\mu$-trace function. \\
   - More generally, for any given $a \in \mathbb{C}^*$, when the fonction $\big(a^k/P_\sigma(a)\big)^\mu$ is well defined  and holomorphic on an open subset $D$  in  $N $,  it is a trace function.\\
     In fact the  function $ \big(\frac{a^k}{P_\sigma(a)}\big)^\mu$ generates in some sense the $\mu$-trace functions and the $\mu$-trace system because  the left ideal $\mathcal{I}_\mu$ may be defined as the annihilator of this function or equivalently  as the annihilator of  each  $\Phi_{\mu, m}, m \in \mathbb{N}$ where $\Phi_{\mu, m}$ is the coefficient of $a^{-m}$ in the Laurent expansion of this function at infinity (in $a$).
      \item In coordinates $z_1, \dots, z_k$ (recall that their elementary symetric functions are $\sigma_1, \dots, \sigma_k$) the partial differential system
      $$(z_i - z_j)\frac{\partial^2}{\partial z_i \partial z_j} -\mu(\partial_{z_i} -\partial_{z_j})$$
for $1 \leq i < j \leq k  $ gives the annihilator of $P_\sigma(z)^{-\mu}:= \big(\prod_{j=1}^k (z - z_j)\big)^{-\mu}$. This corresponds to the singular perturbation
$$ \frac{\partial^2}{\partial z_i\partial z_j} +  \frac{\mu}{z_i - z_j}\Big( \frac{\partial}{\partial z_i} - \frac{\partial}{\partial z_j}\Big) \quad \forall i \not= j \ in \ [1, k]$$
of the partial differential system generated by the operators $ \frac{\partial^2}{\partial z_i\partial z_j}, \forall i\not= j$ which defines trace functions in $D_M^{\frak{S}_k}$, the algebra of $\frak{S}_k$-invariant differential operators on $M := \C^k$ with coordinates $(z_1, \dots, z_k)$.
       \end{enumerate}
       
    \parag{Main results on the $\mu$-Trace System}
    We establish the following results on the regular holonomic  $D_N$-modules $\mathcal{N}_{\mu, \lambda} := D_N \big/( \mathcal{I}_\mu + D_N(U_0 - \lambda))$.
    
    \begin{itemize}
     \item The periodicity Theorem in $\lambda$  Section 3.1.
    \item The minimality Theorem Section 3.2.
    \item The periodicity Theorem in $\mu$ Section 4.2.
    \end{itemize}
    We also describe in section 5  the special cases $\lambda \in \mathbb{Z}$ and $\lambda \in -k\mu + \mathbb{Z}$ when $\mu $ is not in $\mathbb{Z}$.\\
    
    This study  explicites many examples of perverse sheaves associated to representations of the $\pi_1$ of the complement of  the hyper-surface $\{\sigma_k\Delta(\sigma) = 0 \}$ in the affine space with coordinates $\sigma_1, \dots, \sigma_k$, where $\Delta(\sigma)$ is the discriminant of the universal monic  polynomial  of degree $k$, $P_\sigma(z) := z^k + \sum_{h=1}^k (-1)^h\sigma_hz^{k-h}$.\\
    
    This study can be also seen as the study of the first derived direct image with proper support via the projection $pr: N \times \C \to N$ of the $D_{N\times \C}$-module $D_{N\times \C}P_\sigma^{-\mu}$.\\

\section{The $D_N$-modules $\mathcal{W}$ and $\mathcal{M}_\mu$}

\parag{Notations} We define $N := \mathbb{C}^k$ with coordinates $\sigma_1, \dots, \sigma_k$ and we denote by $D_N$ the sheaf on $N$ of holomorphic differential operators on $N$.\\

For basic references on D-modules the reader may consulte \cite{[Bj]}.

\subsection{The $D_N$-module $\mathcal{W}$}

 In this section we recall some basic facts from \cite{[B.24]}. \\
 We shall consider  the $\mathcal{D}_N$-module $\mathcal{W} := \mathcal{D}_N\big/\mathcal{A}$ where $\mathcal{A}$ is the left ideal sheaf in $\mathcal{D}_N$ generated by
\begin{equation}
A_{i, j} := \partial_i\partial_j - \partial_{i+1}\partial_{j-1} \quad {\rm for} \quad i \in [1, k-1] \quad {\rm and} \quad j \in [2, k] 
\end{equation}

Let $\mathcal{D}_N(m)$ be the sub-sheaf of $\mathcal{D}_N$ of partial differential operators of order at most equal to $m$. 
Then let $\mathcal{W}(m)$ be the sub-$\mathcal{O}_N$-module in $\mathcal{W}$ of the classes induced by germs in $\mathcal{D}_N(m)$. \\ 
As we have $\mathcal{A}(m) := \mathcal{D}_N(m) \cap \mathcal{A} = \sum_{i, j} \mathcal{D}_N(m-2)A_{i, j}$ for each $m \in \mathbb{N}$, the quotient $\mathcal{W}(m) = \mathcal{D}_N(m)\big/\mathcal{A}(m) $ injects in $\mathcal{W}$ and we have
$$\mathcal{W} = \cup_{m \geq 0} \mathcal{W}_m.$$
Note that $\mathcal{A}(1) = 0$ \  so $\mathcal{W}(1) = \mathcal{D}_N(1)$. $\hfill \square$\\

It is clear that the characteristic variety of the $\mathcal{D}_N$-module $\mathcal{W}$ is equal to $N \times S(k)$ in the cotangent bundle $ T_N^* \simeq N \times \C^k$ of $N$,  where $S(k)$ is the algebraic cone in $\C^k$ defined by the equations
\begin{equation}
 \eta_i\eta_j - \eta_{i+1}\eta_{j-1} = 0 \quad \forall i \in [1, k-1] \quad {\rm and} \quad  \forall j \in [2, k].
 \end{equation}
This two-dimensional cone $S(k)$ and the corresponding ideal are described in the appendix of \cite{[B.24]}.\\

\parag{Notation} For $\alpha \in \mathbb{N}^k$ we denote $\vert \alpha\vert := \sum_{j=1}^k \alpha_j$ and $w(\alpha) := \sum_{j=1}^k j\alpha_j$. We always has $\vert \alpha\vert \leq w(\alpha) \leq k\vert\alpha\vert$.

\begin{defn}\label{bi-homog.}
Let  $P$ be a germ of section  of $\mathcal{D}_N$. We say that $P$ is {\bf bi-homogeneous of type $(q, r)$} if we may write 
$ P = \sum_{\vert \alpha\vert = q, w(\alpha) = r} a_\alpha\partial^\alpha $ where $a_\alpha$ are  germs of holomorphic functions in $N$.
\end{defn}

It is clear that  any germ $P$  of section of $\mathcal{D}_N$ has a unique decomposition
 $P = \sum_{q, r} P_{q, r}$  where $P_{q, r}$ is a  bi-homogeneous germ of section of $\mathcal{D}_N$  of type $(q, r)$. Note that this sum is finite because for a given order $q$ the corresponding type $(q, r)$ has non zero representative only when $r$ is in $[q, kq]$.\\
 
 The following easy Lemmas are proved in \cite{[B.24]} Lemma 2.1.4 and 2.1.5

\begin{lemma}\label{tres simple}
Let $P$ be a germ of  section of $\mathcal{D}_N$ and write the decomposition of $P$ in  its bi-homogeneous components  as $P = \sum_{q,r}  P_{q, r} $. Then $P$ is a germ of section in $\mathcal{A}$ if and only if for each type $(q, r)$ $P_{q, r}$ is a germ of section in $ \mathcal{A}$.$\hfill \blacksquare$
\end{lemma}

\begin{lemma}\label{base}
 The class  of \ $\partial^\alpha$ in $\mathcal{W}$ only depends on $q := \vert \alpha\vert$ and $r := w(\alpha)$. It will be denoted $y_{q, r} $. Moreover, if $\mathcal{W}_q$ is the sub$-\mathcal{O}_N$-module of $\mathcal{W}$ generated by the $y_{q, r} $ for $r \in [q, kq]$, $\mathcal{W}_q$ is a free $\mathcal{O}_N$-module of rank $kq - q + 1$  with basis  $y_{q, q}, y_{q, q+1}, \dots, y_{q, kq}$ and, as $\mathcal{O}_N$-module, we have the direct decompositions: 
\begin{equation}
 \mathcal{W}(m) = \oplus_{q=0}^m \mathcal{W}_q \quad {\rm and} \quad  \mathcal{W} = \oplus_{q \in \mathbb{N}} \ \mathcal{W}_q .\qquad\qquad  \qquad\qquad \blacksquare \\
 \end{equation}
 \end{lemma}
 
 The global polynomial solutions of  the $\mathcal{D}_N$-module $\mathcal{W}$ are described by Lemma \ref{sol. W} below.

\begin{defn}\label{les mqr} For each $q \in \mathbb{N}$ and each $r \in [q, kq]$ define the polynomial
\begin{equation}
\tilde{m}_{q, r}(\sigma) := q! \sum_{\vert\alpha\vert = q, w(\alpha) = r} \frac{\sigma^\alpha}{\alpha !}
\end{equation}
where $\alpha ! := \prod_{j=1}^k (\alpha_j !)$.
\end{defn}
 Since $\vert \alpha\vert !\big/\alpha !$ is a non negative integer for each $\alpha \in \mathbb{N}^k$, $\tilde{m}_{q, r}(\sigma) $ has integer coefficients: If $\beta := (\alpha, p)$ the equality  $\vert \beta\vert !\big/\beta ! = (\vert \alpha\vert + p)!\big/ \alpha ! p! = C_{\vert \beta\vert}^p \vert \alpha\vert !\big/ \alpha ! $ gives the induction step....\\

The next Lemma is proved in \cite{[B.24]} Lemma 2.1.7.

\begin{lemma}\label{sol. W}
Any $\tilde{m}_{q, r} \in \C[\sigma_1, \dots, \sigma_k]$  is annihilated by the left ideal $\mathcal{A}$ in $\mathcal{D}_N$ and if a polynomial  $P \in \C[\sigma_1, \dots, \sigma_k]$ is annihilated by $\mathcal{A}$, $P$ is, in a unique way, a $\C-$linear combination of the $\tilde{m}_{q, r}$ for $q \geq 0$ and $r \in [q, kq]$ which gives the bi-homogeneous decomposition of $P(\partial_1, \dots, \partial_k) \in \C[\partial_1, \dots, \partial_k]$.$\hfill \blacksquare$
\end{lemma}

Note that an easy way to prove the previous lemma is to remark that the $A_{p, q}$ have pure weight $(p+q)$ and that  $\tilde{m}_{q, r}$ is the sum of the  terms $t$ in the polynomial $(\sigma_1 + \dots + \sigma_k)^q$  such that  $w(t) = r$ and that the polynomial $(\sigma_1 + \dots + \sigma_k)^q$  is homogeneous of degree $q$ and is a solution of $\mathcal{W}$.\\
Note also that $(-1)^r\tilde{m}_{q, r}$ is the sum of the  terms \ $t$ \  such that  $w(t) = r$  in the polynomial $(\sum_{h=1}^k (-1)^h\sigma_h)^q$.

\begin{defn} For each complex number $\mu$ and  integer $q \geq 0$ define $\Phi_{\mu, q}$ as the weight $q$ part of the (convergent) Taylor expansion at $\sigma = 0$ of the holomorphic function
$$ (1 + \sum_{h=1}^k (-1)^h\sigma_h )^{-\mu}.$$
\end{defn}
Then the formula $(1 + x)^{-\mu} := \sum_{q=0}^\infty  \Gamma_{\mu, q}x^q $ with $x := \sum_{h=1}^k (-1)^h\sigma_h$ gives, by the definition of $\Phi_{\mu, q}$, that we have
\begin{equation*}
 \Phi_{\mu, q} =  (-1)^q \sum_p  \Gamma_{\mu, p}\tilde{m}_{p, q} \tag{$\Phi 0$} 
 \end{equation*}
where the sum is finite since for $q$ given $\tilde{m}_{p, q} $ vanishes when $p \geq q+1$.

\parag{Remark} For $\sigma \in N$ and $a \in \mathbb{C}$ let $P_\sigma(a) := a^k + \sum_{h=1}^k (-1)^h\sigma_ha^{k-h}$. Then for $a \in \C^*$, replacing $\sigma_h$ by $\sigma_h a^{-h}$ we find that
$1 + \sum_{h=1}^k (-1)^h \sigma_h$ becomes $P_{\sigma}(a)/a^k$ so that $\Phi_{\mu, q}(\sigma)$ is the coefficient of $a^{-q}$ in the Laurent expansion at $a = 0$ of $P_{\sigma}(a)/a^k$.\\

\subsection{The $\mathcal{D}_N$-module $\mathcal{M}_\mu$}

\begin{defn}\label{ideal}
Let $\mu$ be a complex number and $m \in [2, k]$ be an integer and define the second order differential operator $T^m_\mu$  in the Weyl algebra $\C[\sigma]\langle \partial \rangle$ by
\begin{equation}
\mathcal{T}^m_\mu := \partial_1\partial_{m-1} + \partial_m(E + \mu)  \quad {\rm for} \ m \in [2,k], \  {\rm where} \quad E := \sum_{h=1}^k \sigma_h\partial_h
\end{equation}
Then define the left ideal $\mathcal{I}_\mu$ in $\mathcal{D}_N$ as 
\begin{equation}
\mathcal{I}_\mu := \mathcal{A} + \sum_{m=2}^k \mathcal{D}_N\mathcal{T}^m_\mu
\end{equation}
and let $\mathcal{M}_\mu$ be  the $\mathcal{D}_N$-module 
 \begin{equation}
\mathcal{M}_\mu := \mathcal{D}_N\big/\mathcal{I}_\mu
\end{equation}
\end{defn}

\begin{prop} \label{ZZ2} For $\mu \in \mathbb{C}^*$ the left  ideal $\mathcal{I}_\mu$ generated by $\mathcal{A}$ and the $T^m_\mu$, for $  m \in [2, k]$, is the left  ideal of  germs of partial differential operators   annihilating each $\big(a^k\big/P_\sigma(a)\big)^\mu$ when $\vert a\vert$ is large enough.\\
So $\mathcal{I}_\mu$ is the left ideal in $D_N$ of germs of partial differential operators killing each polynomial  $\Phi_{\mu, q} \in \mathbb{C}[\sigma_1, \dots, \sigma_k]$ for $q \in \mathbb{N}$ which are  defined in formula $(\Phi 0)$ (see also $(\Phi1)$  below).
\end{prop}

\parag{Proof} We first verify that our generators kill $P_\sigma(a)^{-\mu}$ for each $a$, which is the same than to kill $\big(a^k\big/P_{\sigma}(a)\big)^\mu$. \\
Note that it is equivalent, thanks to the Analytic Extension Theorem, to ask that our generators kill $P_\sigma(a)^{-\mu}$ for all $a$ with  $\vert a\vert$ large enough (or even for $a$  in a non discrete subset of the Riemann sphere minus $\{0\}$).\\
Denote by $\mathcal{J}_\mu$ the annihilator of the $P_\sigma(a)^{-\mu}$ for each $a$. Then the formulas\footnote{The first formula below may be taken as a definition of the action of $D_N$ on the expression $ P_{\sigma}(a)^{-\mu}$. Of course, for $\vert a\vert$ large enough compared to $\sigma$ we may define the holomorphic function $Log(a^k/P_\sigma(a))$ and then 
$(a^k/P_{\sigma}(a))^{\mu} := \exp\big(\mu Log(a^k/P_{\sigma}(a)\big)$. Then the natural action of $D_N$ on this holomorphic function (with $a$ constant) is compatible with this formula !}
 \begin{align*}
 & \partial_h(P_\sigma(a)^{-\mu}) = -\mu(-1)^ha^{k-h}P_\sigma(a)^{-\mu-1} \qquad {\rm and \ so} \\
 & \partial_j\partial_h(P_\sigma(a)^{-\mu}) = \mu(\mu+1)(-1)^{j+h}a^{2k-j-h}P_\sigma(a)^{-\mu-2} 
 \end{align*}
 show that, since  $\partial_j\partial_h(P_\sigma(a)^{-\mu}) $   only depends on $j + h$, we have $A_{p, q} \in \mathcal{J}_\mu$ for each pair of integers  $(p, q) \in [1, k-1]\times [2, k]$.\\
 We have also (recall that $E := \sum_{h=1}^k \sigma_h\partial_h $):
 \begin{align*}
 &  E(P_\sigma(a)^{-\mu}) = -\mu\big(\sum_{h=1}^k (-1)^h \sigma_h a^{k-h}\big)P_\sigma(a)^{-\mu -1} \\
 &  E(P_\sigma(a)^{-\mu}) =  -\mu(P_\sigma(a) - a^k)P_\sigma(a)^{-\mu-1}\qquad {\rm and \ so} \\
  & (E+\mu)(P_\sigma(a)^{-\mu}) = \mu a^kP_\sigma(a)^{-\mu-1}.
 \end{align*}
 Then
 $$ \partial_mE(P_\sigma(a)^{-\mu}) = (-1)^m\mu^2a^{2 - m}P_\sigma(a)^{-\mu-1} - (-1)^m\mu(\mu + 1)a^{2k-m}P_\sigma(a)^{-\mu-2}. $$
and we obtain
 \begin{align*}
 &  (T^m + \mu \partial_m)(P_\sigma(a)^{-\mu}) = (-1)^m\mu(\mu+1)a^{2k-m}P_\sigma(a)^{-\mu-2} +  (-1)^m\mu^2a^{k - m}P_\sigma(a)^{-\mu-1} +  \\
 &  \quad  - (-1)^m\mu(\mu + 1)a^{2k-m}P_\sigma(a)^{-\mu-2}   -\mu^2(-1)^ma^{k-m}P_\sigma(a)^{-\mu-1} \qquad {\rm and \ so}\\
 &  (T^m + \mu \partial_m)(P_\sigma(a)^{-\mu}) = 0.
 \end{align*}
 
 So we obtain that $\mathcal{I}_\mu \subset \mathcal{J}_\mu$. The converse is easy since any $\Pi \in \mathcal{J}_\mu$ has its symbol which vanishes on the co-normal space of the hyper-surfaces $\mathcal{H}_a = \{ P_\sigma(a) = 0\}$ and  $Z$ is the union of these co-normal spaces. So the symbol of $\Pi$ vanishes on $Z$.\\
Now,  assuming that $\mathcal{I}_\mu$ is not equal to $\mathcal{J}_\mu$,  choose a $\Pi \in \mathcal{J}_\mu \setminus \mathcal{I}_\mu$ of minimal order. Since its symbol vanishes on $Z$ (which is reduced) there exists $P \in \mathcal{I}_\mu$ with the same symbol. This implies that $\Pi - P $ is in $\mathcal{I}_\mu$ because it has an order strictly less than the order of $\Pi$. Contradiction ! And so $\mathcal{J}_\mu = \mathcal{I}_\mu$.\\

For $\vert a\vert \gg  Sup_{h=1}^k \vert \sigma_h\vert^{1/h}$ the function $Log\big(a^k/P_\sigma(a))$ is holomorphic near $a = \infty$ and so
$(a^k/P_\sigma(a))^\mu := exp(\mu Log(a^k/P_\sigma(a)) $ has a Laurent expansion in $a^{-1}$  at infinity :
\begin{equation*}
 (a^k/P_\sigma(a))^\mu = \sum_{q =0}^\infty  \Phi_{\mu, q}(\sigma)a^{-q}  \tag{$\Phi 1$}
 \end{equation*}
where $\Phi_{\mu, q}$ is a quasi-homogeneous polynomial in $\mathbb{C}[\sigma]$ of weight $q$ (here $\sigma_h$ has weight $h$ for each $h \in [1, k]$). \\
This last point is a consequence of the invariance of $a^k/P_\sigma(a)$ by the $\mathbb{C}^*$-action of $N \times \mathbb{C}$ defined by $\sigma_h \mapsto \rho^h\sigma_h, a \mapsto \rho a$.\\
Then a partial differential operator $\Pi \in D_N$ kills $P_\sigma(a)^{-\mu}$ for any $a$ fixed (large enough) if and only if it kills each $\Phi_{\mu, q}$ for any $q \geq 0$.$\hfill \blacksquare$\\
\parag{Remarks} 
\begin{enumerate}
\item For each $\mu \in \mathbb{C}^*$ and each $q \in \mathbb{N}$ the polynomial $\Phi_{\mu, q}$ is the weight $q$ part of the function $(1+ \sum_{h=1}^k  (-1)^h\sigma_h)^{-\mu}$ (the value for $a = 1$ of  $\big(a^k/P_\sigma(a)\big)^\mu$) which is holomorphic near $\sigma = 0$.\\
This implies that $\Phi_{\mu, q}$ is a linear combination of the polynomials $\tilde{m}_{p, q}$ for $p$ in $[q/k, q]$. The Taylor expansion of the holomorphic function $(1+ \sum_{h=1}^k  (-1)^h\sigma_h)^{-\mu}$ at $\sigma = 0$ gives
\begin{equation*}
(-1)^q q! \Phi_{\mu, q} = \sum_{p \in  [q/k, q]}  (-1)^p \frac{\mu(\mu+1) \dots (\mu+p-1)}{p!} \tilde{m}_{p, q} . \tag{$\Phi 2$}
\end{equation*}

Compare with J. Varouchas formula\footnote{Private communication.}
 \begin{equation}
  N_n = (-1)^n n \sum_{q =1}^n (-1)^q \frac{\tilde{m}_{q, n}}{q} \quad \forall n \in \mathbb{N}^*.
  \end{equation}
  
  \item The following formulas are easily obtained for each $q \in \mathbb{N}$
  \begin{align*}
  & \partial_h\Phi_{\mu, q} = (-1)^{h-1}\mu\Phi_{\mu+1, q-h} \quad \forall h \in [1, k] \quad {\rm with} \quad \Phi_{\mu, q} = 0 \quad {\rm for} \quad q < 0 \\
  & (E + \mu)\Phi_{\mu, q} = \mu \Phi_{\mu+1, q} 
  \end{align*}

\item The formula $(\Phi 2)$ may be obtain directly as follows:
The polynomial $\tilde{m}_{p, q}$ is the weight $q$ part of the polynomial $(\sum_{h=1}^k \sigma_h)^p$ and $\Phi_{\mu, q}$ is the weight $q$ part of
 $(1 + \sum_{h=1}^k (-1)^h\sigma_h)^{-\mu}$. Then the Taylor series of $(1 + x)^{-\mu}$ at $x = 0$ allows to conclude.
\item  The following formulas are easy consequences of the previous remarks:
 $$ \sum_{q+q' = Q} \Phi_{\mu, q}\Phi_{\nu, q'} = \Phi_{\mu+\nu, Q} \quad \forall (\mu, \nu) \in (\mathbb{C} \setminus -\mathbb{N})^2$$
and 
$$ \sum_{q+q' = Q} \tilde{m}_{p, q}\tilde{m}_{p', q'} = \tilde{m}_{p+p', Q} \quad \forall (p, p') \in \mathbb{N}^2.$$
\item For each open set $U$ in $N$ a  section $\Pi$ in  $\Gamma(U, D_N)$ is in $ \Gamma(U, \mathcal{I}_\mu)$ if and only if 
  its annihilates each $(\Phi_{\mu, q})_{\vert U}$ for each $q \in \mathbb{N}$.
\item The characteristic variety of $D_N\big/\mathcal{I}_\mu$ is still the reduced analytic subet $Z$. Its reduced ideal is generated by the symbols of the $A_{p, q}$ and of the $T^m_\mu$.
\item A consequence of the previous proposition is that $\mathcal{M}_\mu$ has no $\mathcal{O}_N$-torsion:\\
 if a differential operator $\Pi$ is such that  $f\Pi$ belongs to $\mathcal{I}_\mu$ for some non zero  $f \in \mathcal{O}_N$ then $\Pi$ kills each $\Phi_{\mu, q}$ for $q \in \mathbb{N}$ on the dense open subset $\{f \not= 0\}$ and so $\Pi$ belongs to $\mathcal{I}_\mu$.$\hfill \square$
\end{enumerate}

  \parag{Notation} For each integer $p = -1, 0, 1$ we note $U_p$ the vector field on $N$ which is the image by the quotient map
  $$ q : \C^k := M \to N := \C^k/\frak{S}_k \simeq \C^k $$
  of the vector field $\sum_{j=1}^k z_j^{p+1}\partial/\partial z_j $. An easy computation\footnote{using, for instance the formulas $\partial \sigma_h/\partial z_j = \sigma_{h-1}(\hat{j}) = \sigma_{h-1} -z_j\sigma_{h-2}(\hat{j})$, where $\sigma_p(\hat{j})$ is the $p$-th symmetric function of $(z_1, \dots, \widehat{z_j}, \dots, z_k)$.} gives,  where $\partial_1, \dots, \partial_k$ denote  the partial derivatives in the coordinates $\sigma_1, \dots, \sigma_k$ of $N$:
  \begin{align*}
  & U_{-1} := \sum_{h=0}^{k-1}  (k-h)\sigma_h\partial_{h+1}, \quad  U_0 := \sum_{j=1}^k h\sigma_h\partial_h, \ {\rm with \ the \ convention} \quad  \sigma_0 = 1, \\
  & U_1 = \sum_{h=1}^k (\sigma_1\sigma_h - (h+1)\sigma_{h+1})\partial_h \quad {\rm with \ the \ convention} \quad \sigma_{k+1} = 0 \quad {\rm and \ so} \\
  & U_1 = \sigma_1E -  \sum_{h=1}^{k-1}  (h+1)\sigma_{h+1}\partial_h.
  \end{align*}
 
 \begin{prop}\label{31/05 suite 1} Assume that $\mu \in \C^*$. Then the left ideal $\mathcal{I}_\mu$ has the following properties:
 \begin{enumerate}
  \item  The characteristic variety of the quotient $\mathcal{M}_\mu : = D_N \big/\mathcal{I}_\mu$ is the reduced and  irreducible  subspace  $Z$ defined by the annulation of the symbols of the  differential operators $A_{p, q}$ and of the $T_\mu^m$.
  \item For each $h \in [1, k]$ the right product by $\partial_h$ sends $\mathcal{I}_{\mu+1}$ in $\mathcal{I}_\mu $.
 \item Each right product by $U_{-1}, U_0$ and  $(U_1+ \mu\sigma_1)$ sends $\mathcal{I}_\mu$ into itself.
 \item The right product by $E + \mu$ sends $\mathcal{I}_{\mu+1}$ into $\mathcal{I}_\mu$.
 \end{enumerate}
 \end{prop}
 
 \parag{Proof} The point $1.$ is  easy, since the symbols of the generators of $\mathcal{I}_\mu$ are the same than the symbols of the generators of the ideal $\mathcal{I}$ introduced in \cite{[B.24]} (which corresponds to $\mu = 0$), and they generate the reduced ideal of  $Z$ which is a pure $(k+1)$-dimensional algebraic sub-set with $1$-dimensional fibers over $N$ (see \cite{[B.24]} Section 3.2). 
  Note that $\mathcal{M}_\mu := D_N\big/\mathcal{I}_\mu$ is sub-holonomic.\\ 
  The point $2.$ is a simple consequence of the fact that $\partial_h$ commutes with $A_{p, q}$ for each $(p, q)$ and the commutation relations
$$ [\partial_mE, \partial_h] = -\partial_h\partial_m \quad {\rm which \ implies} \quad T^m_\mu\partial_h = \partial_hT^m_{\mu-1} \quad \forall h \in [1, k], \forall m \in [2, k].$$
To prove point $3.$ let us compute the action of $U_{-1}, U_0$ and $U_1$ on $P_\sigma(a)$. 
\begin{align*}
& U_{-1}(P_\sigma(a)) = \sum_{h=0}^{k-1} (-1)^{h+1}(k-h)\sigma_h a^{k-h-1} = -P'_\sigma(a) \quad {\rm where} \quad P'_\sigma(a) := \frac{\partial P_\sigma(a)}{\partial a}.\\
& U_0(P_\sigma(a)) = \sum_{h=1}^{k}  (-1)^hh\sigma_ha^{k-h} = -\sum_{h=0}^k (-1)^h(k-h)\sigma_ha^{k-h} + kP_\sigma(a) \\
& U_0(P_\sigma(a)) = -aP'_\sigma(a) + kP_\sigma(a).\\
& E(P_\sigma(a)) = P_\sigma(a) - a^k.
\end{align*}
and so 
\begin{align*}
&  U_1 (P_\sigma(a)) = \sigma_1(P_\sigma(a) - a^k) - \sum_{h=1}^{k-1} (-1)^h(h+1)\sigma_{h+1}a^{k-h} \\
&  U_1 (P_\sigma(a))  = \sigma_1P_\sigma(a) - \sigma_1 a^k +  \sum_{p=0}^k (-1)^pp\sigma_pa^{k-p+1} + \sigma_1a^k \\
&  U_1 (P_\sigma(a))   = \sigma_1P_\sigma(a) + kaP_\sigma(a) -  a^2P'_\sigma(a); 
\end{align*} 
Finally we obtain
\begin{align}
& U_{-1}(P_\sigma(a)^{-\mu}) = \mu P'_\sigma(a)P_\sigma(a)^{-\mu-1}  \\
& U_0(P_\sigma(a)^{-\mu}) =  -k\mu P_\sigma(a)^{-\mu} + \mu aP'_\sigma(a)(P_\sigma(a)^{-\mu-1}) \\
& U_1(P_\sigma(a)^{-\mu}) = -\mu(\sigma_1 + ka)P_\sigma(a)^{-\mu} + \mu a^2P'_\sigma(a)P_\sigma(a)^{-\mu-1}. 
\end{align}
Now, if \, $\Pi(P_\sigma(a)^{-\mu}) = 0$ for $(\sigma, a)$ in a domain $D\subset N \times \mathbb{C}$, we have also for $\mu \not= 0$
 $$\Pi(P'_\sigma(a)P_\sigma(a)^{-\mu-1}) = 0$$
 since $\Pi$ is independent of $a$ it commutes with $\partial/\partial a$ and since
\begin{equation*}
 \frac{\partial}{\partial a} (P_\sigma(a)^{-\mu}) = -\mu P'_\sigma(a)(P_\sigma(a))^{-\mu-1}
 \end{equation*}
we obtain on $D$ the vanishing of $\Pi$ on 
$$U_{-1}(P_\sigma(a)^{-\mu}), \quad  U_0(P_\sigma(a)^{-\mu}) \quad {\rm and \ on} \quad  (U_1 + \mu\sigma_1)(P_\sigma(a)^{-\mu})$$
 proving $3$.\\
 Point $4$ is consequence of the formula $E(P_\sigma(a)) = P_\sigma(a) - a^k$ which gives
 \begin{equation}
 (E + \mu)[P_\sigma(a)^{-\mu}] = \mu a^k(P_\sigma(a))^{-\mu-1}.
 \end{equation}
 and the conclusion follows.$\hfill \blacksquare$\\
 
 \begin{lemma}\label{remplace}
Let $\Pi$ be a non zero germ of  section  of the sheaf $\mathcal{I}_\mu$. Assume that $\Pi$ has order at most $1$. Then $\Pi = 0$.
\end{lemma}

\parag{Proof} Let $\Pi = a_0 + \sum_{h=1}^k a_h\partial_h $. Since $\partial_hP_\sigma(a) = (-1)^ha^{k-h}$ we obtain
$$ \big(a_0P_\sigma(a) - \mu\sum_{h=1}^k (-1)^h a_ha^{k-h}\big)P_\sigma(a)^{-\mu-1} = 0 $$
at least for $\vert a\vert$ large enough. Then we obtain first that $a_0 = 0$ since the polynomial $P_\sigma$ is monic of degree $k$ and 
 then we conclude  that each $a_h$ for $h \in [1, k]$ is zero since $\mu \not= 0$.$\hfill\blacksquare$\\

  For $\mu \in \mathbb{C}^*$ we denote by $H_\mu$ the left ideal in $D_N$ generated by $\partial_1, \dots, \partial_{k-1}$ and $ E + \mu$ where $E := \sum_{h=1}^k \sigma_h\partial_h$. This ideal is the  annihilator of $\sigma_k^{-\mu}$ and we shall denote by $K_\mu$ the quotient $D_N\big/H_\mu$ which is regular holonomic   with characteristic variety the co-normal to the hyper-surface $\{\sigma_k = 0 \}$.\\

 The proof of the  following easy lemma is left to the reader.
   
   \begin{lemma}\label{2/9}
   The left ideal $H_\mu$ in $D_N$ generated by $\partial_1, \dots, \partial_{k-1}, \sigma_k\partial_k + \mu$ where $\mu$ is in $\C \setminus \mathbb{Z}$ is maximal.
   The quotient $K_\mu := D_N/H_\mu$ is a simple holonomic regular $D_N$-module. $\hfill \blacksquare$
   \end{lemma}
   
  \parag{Remark} We have $\mathcal{I}_\mu \subset H_\mu$ and so there is a natural surjective $D_N$-linear map
   $$ \nu_0 : \mathcal{M}_\mu \to K_\mu.$$
   More generally, for each $a \in \C$ the left ideal $H_\mu(a)$  in $D_N$ generated by $\partial_h + a\partial_{h+1}$ for $h \in [1, k-1]$ and $(E + \mu) + a\partial_1$ which is the annihilator
   of $P_\sigma(a)^{-\mu}$,  contains  $I_\mu$ and we have a natural $D_N$-linear surjective map $ \nu_a : \mathcal{M}_\mu \to K_\mu(a) := D_N\big/H_\mu(a)$. For $a = 0$ we have $P_\sigma(0) = (-1)^k\sigma_k$ and since $U_0(\sigma_k) = k\sigma_k$ we find that $\nu_0 : \mathcal{M}_\mu \to K_\mu(0)\simeq D_N\big/(\sum_{h=1}^{k-1} D_N\partial_h + D_N(E + \mu)$ induces a surjective map
   $$ \mathcal{N}_{\mu, -k\mu} := D_N\big/(\mathcal{I}_\mu + D_N(U_0 + k\mu)) \to K_\mu(0).$$
   This map will appear in section 5.

\subsection{Some more solutions of $\mathcal{M}_\mu$}

When $a \to 0$ the quotient $(-1)^k\sigma_k/P_\sigma(a)$ goes uniformly to $1$ assuming that  $\sigma$ stays in a compact set outside the hyper-surface $\{\sigma_k = 0 \}$. This allows to define the logarithm of this quotient on $D \times \{\vert a\vert  < \varepsilon \}$ where $D$ a simply connected domain relatively compact  in $N \setminus \{ \sigma_k = 0 \}$ and where $\varepsilon$ is a real positive number, small enough (depending on $D$). \\
Then, on such a domain we may define, for each $\mu \in \C$ and each $q \in \mathbb{N}$, assuming that $\sigma$ stays in $D$ 
\begin{equation}
\Psi_{\mu, -q}(\sigma) := \frac{1}{2i\pi}\int_{\vert a\vert = \varepsilon/2} a^{-q-1}\Big(\frac{(-1)^k\sigma_k}{P_\sigma(a)}\Big)^\mu da 
\end{equation}
The Cauchy formula shows that $\Psi_{\mu, q}(\sigma)$  is equal to the coefficient of $a^q$ in the Taylor expansion at $a = 0$ of the  function 
$$F_\mu(\sigma, a) := \Big(\frac{(-1)^k\sigma_k}{P_\sigma(a)}\Big)^\mu := \exp\big(-\mu Log(1 + \sum_{h=1}^{k-1} (-1)^{k-h} \frac{\sigma_h}{\sigma_k} a^{k-h})\big). $$ 

So   $\Psi_{\mu, -q}(\sigma) $ is equal to $\frac{1}{q!}\frac{\partial^q F_\mu}{\partial a^q}(\sigma, 0)$ which is a polynomial of degree at most equal to  $q$ in the variables $\sigma_h/\sigma_k$ for $h \in [1, k-1]$. This proves the following lemma, since the function $\Psi_{\mu, -q}$ changes to $\lambda^{-q}\Psi_{\mu, -q}$ by  the change of variables
$$ a \mapsto \lambda a, \quad \sigma_h \mapsto \lambda^h \sigma_h \quad  \forall h \in [1, k]$$
for $\lambda$ near enough to $1$.

\begin{lemma}\label{20/6}
For each $\mu \in \mathbb{C}$, $\sigma_k^q\Psi_{\mu, -q}(\sigma)$ is the restriction to $D$ of a polynomial $Q_{\mu, (k-1)q}$  in $\mathbb{C}[\sigma]$ of pure weight $(k-1)q$. \end{lemma}

\parag{Proof} It is clear that for $\sigma_k \not= 0$ the restriction to $D$ of the function $\Psi_{\mu, -q}(\sigma)$ is a polynomial in $\sigma_1/\sigma_k, \dots, \sigma_{k-1}/\sigma_k$ and that $\sigma_k^q\Psi_{\mu, -q}(\sigma)$ is a polynomial.$\hfill \blacksquare$\\

The interest of these functions is shown by the following result.

\begin{lemma}\label{17/6 1}
Define on  a domain $D$ as above the function
\begin{equation}
\tilde{\Phi}_{\mu, -q}(\sigma) := \big((-1)^k\sigma_k\big)^{-\mu}\Psi_{\mu, -q}(\sigma) =  \big((-1)^k\sigma_k\big)^{-\mu}\sigma_k^{-q}Q_{\mu, (k-1)q}(\sigma)
\end{equation}
Then $\tilde{\Phi}_{\mu, -q}$ is solution on the simply connected  domain $D \subset\subset N \setminus \{\sigma_k = 0 \}$ of the $D_N$-module $\mathcal{M}_\mu$, where the polynomial $Q_{\mu, (k-1)q} \in \C[\sigma]$ has pure weight $(k-1)q$.  Moreover it satisfies 
\begin{equation}
(U_0 + k\mu + q)(\tilde{\Phi}_{\mu, -q}) = 0
\end{equation}
So In fact  $\tilde{\Phi}_{\mu, -q}$   is a solution on $D$ of the $D_N$-module
 $$\mathcal{N}_{\mu, -k\mu -q} := D_N\big/(\mathcal{I}_\mu + D_N(U_0 + k\mu + q))$$
 which is considered below !
\end{lemma}

\parag{Proof} Since each generator of the left ideal $\mathcal{I}_\mu$ kills $P_\sigma(a)^{-\mu}$ for any $a \in \mathbb{C}$ it is clear that it kills $\big((-1)^k\sigma_k\big)^{-\mu}\Psi_{\mu, -q}(\sigma) $. So $\tilde{\Phi}_{\mu, -q}(\sigma) $ is a solution of $\mathcal{M}_\mu$ on $D$.\\
The formula $U_0(\Psi_{\mu, -q}) = -q\Psi_{\mu, -q}$ is easily obtained (as already explained above)  by the change of variable $a \to \lambda a$ in the integral  and the changes $\sigma_h \to \lambda^h\sigma_h, h \in [1, k]$ which correspond to the action of $U_0$ on $N$. $\hfill \blacksquare$

\parag{Remark} The previous $\mu$-trace functions found above are the analog of the trace functions defined  by the integrals
$$ \Phi_{-q}(\sigma) = \frac{1}{2i\pi} \int_{\vert a \vert = \varepsilon/2} a^{-q} Log\big((-1)^k\sigma_k/P(\sigma, a)\big)\frac{da}{a} $$
After an integration by parts for $q \geq 1$ we obtain
$$ \Phi_{-q}(\sigma) = \frac{1}{2i\pi} \int_{\vert a \vert = \varepsilon/2} a^{-q}\frac{P'(\sigma, a)}{qP(\sigma, a)}da = N_{-q}(\sigma). $$
where $N_{-q}(\sigma) = \sum_{j=1}^k z_j^{-q}$. They are rational functions defined on $N \setminus \{\sigma_k = 0\}$.
More precisely, $\sigma_k^q\Phi_{-q}(\sigma)$ is a polynomial in $\mathbb{C}[\sigma]$ with pure weight $(k-1)q$.\\

\begin{lemma}\label{5/2} Let $D$ be a a simply connected domain relatively compact  in the open set  $N \setminus \{ \sigma_k = 0 \}$ and  let  $\varepsilon$ be a real positive number, small enough (depending on $D$). Then we have on $D$  for each $q \in \mathbb{N}$
\begin{equation*}
 U_{-1}[\tilde{\Phi}_{\mu, -q}] = -(q+1)\tilde{\Phi}_{\mu, -q-1} \quad {\rm and} \quad (E + \mu)(\tilde{\Phi}_{\mu, -q}) = \mu \tilde{\Phi}_{\mu+1, -q+k} \tag{R}
 \end{equation*}
\end{lemma}

\parag{Proof} Using formulas $(16), (20)$ and $(21)$ we obtain
$$  U_{-1}[\tilde{\Phi}_{\mu, -q}] (\sigma) = \frac{\mu}{2i\pi} \int_{\vert a \vert = \varepsilon/2} a^{-q-1}\frac{P'(\sigma, a)}{P_\sigma(a)^{\mu+1}}da $$
which gives, after integration by parts, the announced formula.
The second formula is a direct consequence of $(19 )$ with an integration by parts.$\hfill \blacksquare$\\

The first formula above  is the analog of the formula 
$$U_{-1}[N_{-q}/q] = -(q+1)N_{-q-1}/(q+1), \quad  \forall q \in \mathbb{N}^*$$
for the trace functions. It allows to compute $\tilde{\Phi}_{\mu, -q}$ as
 $$U_{-1}^q[\tilde{\Phi}_{\mu, 0}] = (-1)^{q}q!\tilde{\Phi}_{\mu, -q}$$
 where  
 $\tilde{\Phi}_{\mu, 0} = ((-1)^k\sigma_k)^{-\mu}$.\\
 The second formula is the analog of the formula
 $$ E(N_q) = qM_q $$
 where $M_q := \Phi_{1, q} := \frac{1}{2i\pi}\int_{\vert a\vert = R} \frac{a^{q+k -1} da}{P_\sigma(a)} $, see \cite{[B.22]} paragraph 5.2.\\
 
 The relation $(22)$ shows that on a simply connected domain  $D \subset \subset N \setminus \{\sigma_k = 0\}$ the function $\tilde{\Phi}_{\mu, -q}$ is solution of the $D_N$-module $\mathcal{N}_{\mu, \lambda} := D_N\big(\mathcal{I}_\mu + D_N(U_0 - \lambda)\big)$  for $\lambda = -k\mu - q$. \\
 
 The relations  $(R)$ correspond to $D_N$-linear maps
 $$ \square U_{-1} : \mathcal{N}_{\mu, \lambda} \to  \mathcal{N}_{\mu, \lambda +1} \quad {\rm and} \quad \square (E + \mu) : \mathcal{N}_{\mu +1, \lambda} \to \mathcal{N}_{\mu, \lambda}  $$
 which are studied below. Note that $-k(\mu+1) - q + k  = -k\mu - q $ so  that the right product by $E + \mu$ preserves the value of $\lambda$ since $E$ and $U_0$ commute.
 
\subsection{Other examples}
 
 Let $\theta : \tilde{N} \to N'$ be the universal cover of $N' := \{ \sigma \in N \ / \ \sigma_K\Delta(\sigma) \not= 0 \}$, and define the hypersurface $\mathcal{H}$ in $\tilde{N} \times \C$ by
 $$ \mathcal{H} := \{(\tilde{\sigma}, a) \in \tilde{N} \times \C \ / \ P_{\theta(\tilde{\sigma})}(a) = 0 \}. $$
 Then $\theta_{\vert \mathcal{H}} : \mathcal{H} \to \tilde{N}$ is a $k$-sheeted (non ramified) cover of the simply connected manifold $\tilde{N}$. So it is the union of $k$ complex manifolds $\mathcal{H}_1, \dots, \mathcal{H}_k$ isomorphic to $\tilde{N}$ via $\theta$. Then $H_1(\tilde{N}\times \C \setminus \mathcal{H}, \mathbb{Z})$  is isomorphic to $\mathbb{Z}^k$, and for each $\gamma $ in this group , for each $q \in \mathbb{N}$ and each $\mu \in \C$ the function on 
 $\tilde{N}$ defined by
 \begin{equation*}
  \Phi_{\gamma, q, \mu}(\tilde{\sigma}) := \int_\gamma a^{q-1}\big (\frac{a^k}{P_{\theta(\tilde{\sigma})}(a)}\big)^\mu da \tag{@} 
  \end{equation*}
 induces a solution of the $\mu$-Trace system on any simply connected open subset $U$ in $N'$, identifying $U$ with the image of a holomorphic section $s : U \to N'$ of $\theta$.\\
 
 Consider now the family parametrized by $N$ of curves in $\mathbb{P}_2(\C)$ which are the compactifications of the affine curves
 $$ C_{p, k, \sigma} := \{ (x, y) \in \C^2 \ / \  y^p x^k = P_{\sigma}(x) \}.$$
 Then, for $\mu = -1/p$ the differential form $x^{q-1}ydx$ on $C_{p, k, \sigma} $ have periods which are given by formulas  analogous to  $(@)$ with $\mu = -1/p$, for $\sigma \in N'$.\\

 \section{The $\mathcal{D}_N$-modules $\mathcal{N}_{\mu, \lambda}$}

 \subsection{Periodicity in $\lambda$}
 
 Define the left ideal $\mathcal{J}_{\mu, \lambda}$ in $D_N$ as the sum $\mathcal{I}_\mu + D_N(U_0 - \lambda)$ for any $(\mu, \lambda) \in \C^2$, and put $\mathcal{N}_{\mu, \lambda} := D_N\big/\mathcal{J}_{\mu, \lambda}$.\\
 Define the left  $D_N$-linear $\mathscr{T}_{\mu, \lambda} : \mathcal{N}_{\mu, \lambda} \to \mathcal{N}_{\mu, \lambda+1}$ by the right product by $U_{-1}$  in $D_N$  and $\mathscr{H}_{\mu, \lambda} : \mathcal{M}_{\mu} \to \mathcal{M}_{\mu}$ by the right product by $U_0 - \lambda$ in $D_N$ for  each $\lambda \in \mathbb{C}$ which are studied below.

 \begin{prop}\label{fond.1}
Assume that $\mu$ is in $\mathbb{C}$. For each $\lambda \in \C$ we have an exact sequence of left $\mathcal{D}_N$-modules on $N$
\begin{equation}
0 \to \mathcal{M}_\mu \overset{\mathscr{H}_{\mu, \lambda}}{\longrightarrow} \mathcal{M}_\mu \overset{q_{\mu, \lambda}}{\longrightarrow} \mathcal{N}_{\mu, \lambda} \to 0 
\end{equation}
where $q_{\mu, \lambda}$ is the obvious quotient map.
\end{prop}

\parag{Proof}  The  quotient map $q_{\mu, \lambda}$ is surjective by definition, so the point is to prove that the kernel of $q_{\mu, \lambda}$ is isomorphic to
 $\mathcal{M}_\mu$.\\ 
This kernel is obviously given by
\begin{equation}
 D_{N}(U_0 - \lambda)\big/\mathcal{I}_\mu\cap D_{N}(U_{0}- \lambda) \overset{i}{\to} \big(\mathcal{I} _\mu + D_{N}(U_0 - \lambda)\big)\big/\mathcal{I} _\mu  \simeq  \mathcal{J}_{\mu, \lambda}\big/ \mathcal{I}_\mu 
\end{equation}
where the fact that $i$ is an isomorphism  is given by  the following lemma.

\begin{lemma}\label{fond.0}
Let $P$ be a germ in $D_{N, \sigma}$ for some $\sigma \in N$ such that $P(U_0 - \lambda) $ is in 
$\mathcal{I}_{\mu, \sigma}$. Then $P$ is in $\mathcal{I}_{\mu, \sigma}$. So \  $\mathcal{I}_{\mu, \sigma} \cap D_{N, \sigma}(U_0 -\lambda) = \mathcal{I}_{\mu, \sigma}(U_0 -\lambda)$.
\end{lemma}

\parag{Proof} Assume that the lemma is wrong. Then let $P_{0}$  in $ D_{N, \sigma}$ having minimal order among germs $P$ in $ D_{N, \sigma}$ satisfying the following properties
\begin{enumerate}
\item $P(U_0 -\lambda)$ is in $\mathcal{I}_{\mu, \sigma}\cap D_{N, \sigma}(U_{0}- \lambda)$ ;
\item $P$ is not in $ \mathcal{I}_{\mu, \sigma}$.
\end{enumerate}
 Let $\pi$ be the symbol of $P_{0}$  and let $g$ be the symbol of $U_{0}$. Then $\pi g $ vanishes on $p^{-1}(V) \cap Z$ where $V$ is a neighborhood of $\sigma$  in $N$ and $p : {T_N}^* \to N$ the projection. But we know that $g$ does not vanish on any non empty open set of $Z$ because $\{g = 0\} \cap Z$ has pure co-dimension $1$ in $Z$ (see \cite{[B.24]} Section 3.2 ). Then $\pi$ vanishes on $(V \times \C^k) \cap Z$ where $V$ is a neighborhood  of $\sigma$ in $N$ and, as we have proved that $Z$ is reduced and is the characteristic cycle of $\mathcal{M}_\mu$, their exists a germ $P_{1}$ in $\mathcal{I}_{\mu, \sigma} $ with symbol  equal to $\pi$. Then $(P_{0}- P_{1})(U_{0}- \lambda)$ satisfies again the properties 1. and  2. and is of order strictly less than the order of $P_{0}$. So  $P_0 - P_{1}$ is in $\mathcal{I}_{\mu, \sigma}$ and this contradicts the fact that we assumed that $P_{0}$ is not in $\mathcal{I}_{\mu, \sigma}$.$\hfill \blacksquare$

\parag{End of proof of \ref{fond.1}} The previous lemma shows that for each $\lambda \in \C$
$$ \mathcal{I}_\mu \cap D_{N}(U_{0} - \lambda) = \mathcal{I}_\mu(U_{0} - \lambda).$$

So the right multiplication by $U_{0}- \lambda$ induces an isomorphism of left $\mathcal{D}_N$-modules
$$ \mathcal{M}_\mu \to D_{N}(U_{0}- \lambda)\big/\mathcal{I}_\mu(U_{0}- \lambda) $$
and the kernel of $q_{\mu, \lambda}$ is isomorphic to $\mathcal{M}_\mu$ by the inverse of this isomorphism.$\hfill \blacksquare$\\

\begin{defn}\label{etoiles}
Define the $\mathcal{D}_N$-module $\mathfrak{N}_\mu$ as the quotient $\mathcal{D}_N\big/(\mathcal{I}_\mu + \mathcal{D}_NU_{-1})$. 
 For each $\lambda \in \C$ then define\footnote{We use here point 3 in Proposition \ref{31/05 suite 1}.} $\mathscr{T}_{\mu, \lambda} : \mathcal{N}_{\mu, \lambda} \to \mathcal{N}_{\mu, \lambda+1}$ as the $\mathcal{D}_N$-linear map induced by $\mathscr{T}_\mu$.
\end{defn}

\begin{lemma}\label{11/06}
For each $\lambda $ and $\mu$ in $\C$ there is a natural  $D_N$-linear isomorphism from the  co-kernel of the $\mathcal{D}_N$-linear map $\mathscr{T}_{\mu, \lambda}$
to  the co-kernel of  the $\mathcal{D}_N$-linear  map 
 $\tilde{\mathscr{H}}_{\lambda+1}  : \mathfrak{N}_\mu \to \mathfrak{N}_\mu $ induced by 
$\mathscr{H}_{\mu,\lambda+1} $.\\
There is also  a natural surjective  map of $\mathcal{D}_N$-modules between the kernels of $\mathscr{T}_{\mu, \lambda}$ and the kernel of $ \tilde{\mathscr{H}}_{\mu, \lambda+1}$.
\end{lemma}

\bigskip

\bigskip

\parag{proof} Consider the commutative diagram of left $\mathcal{D}_N$-modules with exact lines and columns:
$$ \xymatrix{ \quad & \quad & \quad & \quad & 0 \ar[d] & \quad \\ \quad & \quad & \quad & \quad &Ker(\tilde{\mathscr{H}}_{\mu, \lambda+1}) \ar[d] & \quad \\
 \quad & \quad & \mathcal{M}_\mu \ar[d]^{\mathscr{H}_{\mu, \lambda}} \ar[r]^{\mathscr{T}_\mu} & \mathcal{M}_\mu \ar[d]^{\mathscr{H}_{\mu, \lambda+1}} \ar[r] & \mathfrak{N}_\mu \ar[d]^{\tilde{\mathscr{H}}_{\mu, \lambda+1}} \ar[r] & 0 \\ \quad & \quad & \mathcal{M}_\mu \ar[r]^{\mathscr{T}_\mu}  \ar[d] & \mathcal{M}_\mu \ar[d] \ar[r] & \mathfrak{N}_\mu \ar[r] \ar[d]^\xi & 0 \\ 
0 \ar[r] & Ker(\mathscr{T}_{\mu, \lambda}) \ar[r]  \ar@/^5pc/[rrruuu]^\eta  & \mathcal{N}_{\mu, \lambda} \ar[r]^{\mathscr{T}_{\mu, \lambda}} \ar[d] & \mathcal{N}_{\mu, \lambda+1} \ar[d] \ar[r] &  \mathcal{N}^\square_{\mu, \lambda+1} \ar[r] & 0 \\ \quad & \quad & 0 & 0 & \quad & \quad } $$
where $\mathcal{N}^\square_{\mu, \lambda+1}$ is, by definition, the co-kernel of $\mathscr{T}_{\mu, \lambda} : \mathcal{N}_{\mu, \lambda} \to \mathcal{N}_{\mu, \lambda+1}$ and where the map $\xi$ is defined thanks to the commutation relation $(U_0 -\lambda)U_{-1} = U_{-1}(U_0 - (\lambda + 1))$ which gives also the commutativity of the left   rectangles.\\
 By a  simple diagram chasing it is easy to see that $\mathcal{N}^\square_{\mu, \lambda+1}$  is  also the co-kernel of $\tilde{\mathscr{H}}_{\mu, \lambda+1} :  \mathfrak{N}_\mu \to \mathfrak{N}_\mu$,
and an  elementary diagram chasing gives the existence of  a $D_N$-linear  map $\eta$  between the  kernels of $\mathscr{T}_{\mu, \lambda}$. and  $ \tilde{\mathscr{H}}_{\mu, \lambda+1}$ and the fact that this map is surjective.  $\hfill \blacksquare$\\

We shall prove now that for $\lambda \not=  -1, -k\mu$ the map $\mathscr{T}_{\mu, \lambda}$ is an isomorphism of left $\mathcal{D}_N$-modules. This implies $\mathcal{N}^\square_{\mu, \lambda} = \{0\}$ for $\lambda \not= 0$ and $\lambda \not= 1- k\mu$. 

\begin{thm}\label{isom 2} {\bf [Periodicity in $\lambda$]} For any $\lambda$ and $\mu$ in $\C$ let $\mathscr{G}_{\mu, \lambda} : \mathcal{N}_{\mu, \lambda+1} \to \mathcal{N}_{\mu, \lambda}$ and $\mathscr{T}_{\mu, \lambda} : \mathcal{N}_{\mu, \lambda} \to \mathcal{N}_{\mu, \lambda+1}$ be respectively the left $\mathcal{D}_N$-linear map given by right multiplication by $U_1 + \mu\sigma_1$ and by $U_{-1}$. Then we have   \begin{align*}
  & \mathscr{T}_{\mu, \lambda }\circ \mathscr{G}_{\mu, \lambda} =   \square (U_1+ \mu\sigma_1)U_{-1} =  (\lambda + k\mu)(\lambda + 1) \quad {\rm and} \\
   & \mathscr{G}_{\mu, \lambda}\circ \mathscr{T}_{\mu, \lambda} =  \square U_{-1}(U_1+ \mu\sigma_1) = (\lambda + k\mu)(\lambda + 1) \quad \tag{$T$}
   \end{align*}
on $\mathcal{N}_{\mu, \lambda + 1}$ and $\mathcal{N}_{\mu, \lambda}$ respectively.\\
So for $\lambda \not= -1, -k\mu$ the  left $\mathcal{D}_N$-linear map $\mathscr{T}_{\mu, \lambda}$ is an isomorphism between $\mathcal{N}_{\mu, \lambda}$ and $\mathcal{N}_{\mu, \lambda + 1}$  and $\mathscr{G}_{\mu, \lambda}$ and isomorphism from $\mathcal{N}_{\mu, \lambda+1}$ to $\mathcal{N}_{\mu, \lambda}$.
\end{thm}

\parag{Remark} Note that this result gives also
\begin{itemize}
\item $ Im(\mathscr{G}_{\mu, -1}) \subset Ker(\mathscr{T}_{\mu, -1})$ in $\mathcal{N}_{\mu, -1}$;
\item $Im(\mathscr{T}_{\mu, -1}) \subset  Ker(\mathscr{G}_{\mu, -1})$ in $\mathcal{N}_{\mu, 0}$;
\item  $ Im(\mathscr{G}_{\mu, -k\mu}) \subset Ker(\mathscr{T}_{\mu, -k\mu})$ in $\mathcal{N}_{\mu, -k\mu}$;
\item $Im(\mathscr{T}_{\mu, -k\mu}) \subset  Ker(\mathscr{G}_{\mu, -k\mu})$ in $\mathcal{N}_{\mu, 1 - k\mu}$.
\end{itemize}

\parag{Proof}   Thanks to the formulas $(8), (9), (10)$ and (8-$bis$)  we obtain:
  \begin{align*}
& U_{-1}(\Phi_{\mu, q}) = (k\mu + q - 1)\Phi_{\mu, q-1} \\
& (U_1 + \mu\sigma_1)(\Phi_{\mu, q}) = (q+1)\Phi_{\mu, q+1} \\
& U_0(\Phi_{\mu, q}) = q\Phi_{\mu, q} 
\end{align*}
This iomplies,  using now Proposition  \ref{ZZ2} :
 
  $$ (U_1+ \mu\sigma_1)U_{-1} = U_0(U_0 + k\mu -1)  \quad  {\rm modulo} \  \mathcal{I}_\mu $$
  
  and
  
  $$U_{-1}(U_1+ \mu\sigma_1) = (U_0 + k\mu )(U_0 + 1)   \quad  {\rm modulo} \  \mathcal{I}_\mu $$
  
 This gives,  respectively on  $\mathcal{N}_{\mu, \lambda + 1}$ and on  $\mathcal{N}_{\mu, \lambda}$:
  \begin{align*}
  &  \square (U_1+ \mu\sigma_1)U_{-1} = (\lambda +1)(\lambda + k\mu) \quad {\rm and} \\
   & \square U_{-1}(U_1+ \mu\sigma_1) = (\lambda + k\mu)(\lambda + 1).
   \end{align*}
   concluding the proof. $\hfill \blacksquare$\\
   
   This theorem has the following consequences, assuming\footnote{For $\mu \in \mathbb{N}$ see \cite{[B.24]}} that $\mu \not\in \mathbb{Z}$:
   \begin{enumerate}
   \item For  any $\lambda \not= -1$ and $\lambda \not=  -k\mu $ the $D_N$-modules $\mathcal{N}_{\mu, \lambda+1}$ and $\mathcal{N}_{\mu, \lambda}$ are isomorphic.
    \item For $\lambda = -k\mu +p$ and $p \in \mathbb{N^*}$ the $D_N$-module $\mathcal{N}_{\mu, -k\mu + p}$ is isomorphic to $\mathcal{N}_{\mu, 1 -k\mu}$.
   \item   For $\lambda = -k\mu - q$ and $q \in \mathbb{N}$ the $D_N$-module $\mathcal{N}_{\mu, -k\mu-q}$ is isomorphic to $\mathcal{N}_{\mu, -k\mu}$.
   \end{enumerate}
   
    The structures of the $D_N$-modules  $\mathcal{N}_{\mu, 0},  \mathcal{N}_{\mu, -1}$ and $\mathcal{N}_{\mu, 1-k\mu}, \mathcal{N}_{\mu, -k\mu}$ are examined in Section 5.

  \subsection{Minimality}

 First we recall two results from \cite{[B.24]} (see Corollaries 3.2.2 and 3.2.3 in {\it loc. cit.}) which apply directly to our situation since the symbols of the generators of the ideal $\mathcal{I}_\mu$ are the same that the symbols for the ideal $\mathcal{I}$ in \cite{[B.22]} or  \cite{[B.24]} ( corresponding to $\mu = 0$).
 
 \begin{cor}\label{charact. 00} For any choices of complex numbers  $\mu$ and $\lambda$ the characteristic cycle of $\mathcal{N}_{\mu, \lambda}$ is the cycle associated to the ideal $I_Z + (g)$ in $\mathcal{O}_N[\eta]$ where $g$ is the symbol of $U_0$.\\
 Also the characteristic cycle of  $\mathfrak{N}_\mu$ is the cycle associated to the ideal $I_Z + (\gamma)$ in $\mathcal{O}_N[\eta]$ where $\gamma$ is the symbol of $U_{-1}$.$\hfill \blacksquare$
\end{cor}

\begin{cor}\label{our sheaf case} Let $\lambda$ and $\mu$ be complex numbers and let  $U := U_0 -\lambda$. Then for any non zero germ $Q \in \mathcal{I}_\mu + \mathcal{D}_NU$ of order $q+1$ there exist a germ $P \in \mathcal{I}_\mu$ of order at most $q+1$ and a germ $B \in \mathcal{D}_N$ of order at most $q$ such that $Q = P + BU$. $\hfill \blacksquare$
\end{cor}

\parag{Remarks} 
 \begin{enumerate}
\item  Define for each $q \geq 0$
$$ \mathcal{J}_{\mu, \lambda}(q+1) = \mathcal{I}_\mu(q+1) + \mathcal{D}_N(q)(U_0 - \lambda) .$$
Then the second  corollary above  gives, for each $\lambda \in \C$ and for each $q \in \mathbb{N}^*$ the equality $\mathcal{J}_{\mu, \lambda} \cap \mathcal{D}_N(q) = \mathcal{J}_{\mu, \lambda}(q)$. This implies that the natural map
\begin{equation}
\mathcal{N}_{\mu, \lambda}(q) \to \mathcal{N}_{\mu, \lambda} 
\end{equation}
is injective 
\item Note that $\mathcal{J}_{\mu, \lambda}(0) := \mathcal{I}_\mu(0) = \{0\}$ as no non zero differential operator of order $0$ annihilates the polynomial $\Phi_{\mu, 0} \equiv 1$.
\item Also we have $\mathcal{I}_\mu(1) = \{0 \}$,  thanks to Lemma \ref{remplace},  and  so
 $$\mathcal{J}_{\mu, \lambda}(1) = \mathcal{O}_N(U_0 - \lambda).$$
\end{enumerate}

Note that the computation of the characteristic variety of $\mathcal{N}_{\mu, \lambda}$  is the same than the characteristic variety of $\mathcal{N}_{0, \lambda}$  computed in \cite{[B.24]}. This  shows that the following Corollary (the $\mu$-analog of Corollary  2.3.2  in {\it loc. cit.}\footnote{Proposition 2.3.1  in {\it loc. cit.} is not correct since for any $a \in \C$ the co-normal of the hyper-surface $H_a := \{\sigma \in N \ / \ P_\sigma(a) = 0\}$ is contained in $Z$. But  Corollary 2.3.2 is true thanks to the identification of the characteristic variety of $\mathcal{N}_\lambda$ given later on in Section 3.2  of this article. So this (wrong) proposition may be deleted from \cite{[B.24]} without changing the results of \cite{[B.24]}.}) holds true.

 \begin{cor}\label{co-torsion 2}
 Let $Q$ be a coherent $D_N$-module quotient of $\mathcal{N}_{\mu, \lambda}$ which is supported in a closed analytic subset $S$ in $N$ with empty interior in $N$. If $Q$ vanishes near the generic points of $\{\sigma_k = 0 \} \cup \{\Delta = 0 \}$, then $Q = \{0\}$.$\hfill \blacksquare$\\
 \end{cor}
 
The following adaptation of Lemma 2.2.6 in \cite{[B.24]} is immediate.

\begin{lemma}\label{computation}
Let $q \geq 2$ be an integer, $\alpha \in \mathbb{N}^k$ such that $\vert \alpha \vert = q-2$  and let $m$ be an integer in the interval $[2, k]$. The class induced by $\partial^\alpha T_\mu^m$ in $\mathcal{W}$ only depends on the integers $q$ and $r := w(\alpha) + m$. This class is given by the formula (with the convention $\sigma_0 \equiv 1$)
\begin{equation}
[\partial^\alpha T_\mu^m] = \sum_{h=0}^k \sigma_hy_{q, r+h} + (q + \mu -1)y_{q-1, r},
\end{equation}
where $y_{q, r}$ is the class induced by $\partial^\gamma$ in $\mathcal{W}$ for any $\gamma \in \mathbb{N}^k$ such that $\vert \gamma \vert = q$ and $w(\gamma) = r$ (see Lemma \ref{base}).\\
Let  $\lambda$ be a complex number and  let $\beta \in \mathbb{N}^k$ such that $\vert \beta\vert = q-1$ and $w(\beta) = r$. The class induced by $\partial^\beta(U_0 - \lambda)$ in $\mathcal{W}$ only depends on $\lambda$ and on the integers $q$ and $r$. This class is given by
\begin{equation}
[\partial^\beta(U_0 - \lambda)] = \sum_{h=1}^k  h\sigma_hy_{q, r+h} + (r - \lambda)y_{q-1, r} .
\end{equation}
Also the class induced by $\partial^\beta U_{-1}$ in $\mathcal{W}$, again for $\vert\beta\vert = q-1$ and $w(\beta) = r$,  only depends on the integers $q$ and $r$. This class is given by
\begin{equation}
[\partial^\beta U_{-1} ]= \sum_{h=0}^k (k-h)\sigma_hy_{q, r+h+1}  + (k(q-1) - r)y_{q-1, r +1}
\end{equation}
where, for $r = k(q-1)$, the last term in $(28)$ is equal to $0$ by convention.$\hfill \blacksquare$\\
\end{lemma}

 \begin{defn}\label{5/9}
For $\mu \not= 0$ and $\lambda \in \C\setminus \mathbb{N}^*$ define  for any integer  $q \geq 2$ and for any integer $r \in [q, k(q-1)]$  the following sections of  $\mathcal{W}_q$:
 \begin{align*}
& \theta_{q, r} := (r - \lambda)[\partial^\alpha\mathcal{T}_\mu^m]  - (q+ \mu - 1)[\partial^\beta(U_0 - \lambda)] \tag{$\theta(\mu)$}\\
&{\rm so} \quad   \theta_{q, r} = \sum_{h=0}^k  (r - \lambda - (q+\mu-1)h)\sigma_hy_{q, r+h}. 
\end{align*}
where in Formula $(\theta(\mu))$ we assume that  $\alpha \in \mathbb{N}^k$ and $m \in [2, k]$  satisfy $\vert \alpha \vert = q-2$ and $w(\alpha) = r-m$, and that $\beta \in \mathbb{N}^k$ satisfies $\vert \beta \vert = q-1$ and $w(\beta) = r$.\\

\end{defn}

\begin{cor}\label{charact. 2}
 For any integer $q \geq 1$ the kernel of the quotient map
  $$ l_q : \mathcal{W}(q) \to \mathcal{N}_{\mu, \lambda}(q)$$
  is equal to the sub-$\mathcal{O}_N$-module generated by $U_0 - \lambda \in \mathcal{W}(1)$ and the elements $\theta_{p, r}$, for each $p \in [2, q]$, and for each $ r \in [p, (k-1)p],$. 
\end{cor}

\parag{Proof} We have to prove that if a non zero differential operator $P$ of order $p \leq q$ is in $\mathcal{J}_{\mu, \lambda}$ then it may be written as $Q + B(U_0 -\lambda)$ with $Q \in \mathcal{I}_\mu$ of order at most $p$ (or $Q = 0$) and $B$ of order at most $p-1$ (or $B = 0$). When $p \geq 2$ this is precisely the statement proved in Corollary \ref{our sheaf case}.\\
For $p \leq 1$ the only $P$ which are in $\mathcal{J}_{\mu, \lambda}(1)$ are in 
$\mathcal{O}_N(U_0 - \lambda)$ thanks to  Remark 3 above. $\hfill \blacksquare$\\ 

\begin{lemma}\label{basis}
Let  $\lambda$ be in  $ \C \setminus \mathbb{N}^*$;  for each integer $q \geq 2$ the elements $\theta_{q, r}$ and $y_{q, s}$, with $r \in [q, k(q-1)]$ and $ s \in [k(q-1)+1, kq]$ form a $\mathcal{O}_N-$basis of $\mathcal{W}_q$.
\end{lemma}

\parag{proof} Let $\mathcal{W}_{q, p}$ be the $\mathcal{O}_N$-module of $\mathcal{W}_q$ with basis the $y_{q, r}$ for $r \geq p+1$. Then we have for $r \in [q, k(q-1)]$
$$ \theta_{q, r} \in  (r- \lambda)y_{q, r} + \mathcal{W}_{q, r+1} $$
so the determinant of the $k(q-1) - q + 1 + k  = kq - q + 1$ vectors $\theta_{q, r}, y_{q, s}$ in the basis $\big(y_{q, r}, r \in [q, kq]\big)$ of $ \mathcal{W}_q $ is  upper triangular and is equal to $\prod_{r = q}^{k(q-1)} (r - \lambda) $ which is in $\C^*$ as soon as $\lambda$ is not in the subset $[q, k(q-1)]$ of $\mathbb{N}^*$.$\hfill \blacksquare$\\

The proposition below is the $\mu$-analog of Proposition 3.3.5 in \cite{[B.24]} and the proof is the same using formulas $(\theta(\mu))$ instead of $(\theta(0))$.

\begin{prop}\label{non torsion}
Let  $q \geq 1$ be an integer and assume that $\lambda $ is not an integer in $[0, k(q-1)]$. Let $L_q : \mathcal{W}_q \to \mathcal{N}_{\mu, \lambda}(q)$ be the restriction to $\mathcal{W}_q$ of  quotient map $l_q$. This 
$\mathcal{O}_N$-linear map is surjective and its kernel  is the sub-module of $\mathcal{W}_q$ with basis the $\theta_{q, r}$ for  $r \in [q, k(q-1)]$. So $\mathcal{N}_{\mu, \lambda}(q)$ is a free $\mathcal{O}_N$-module of  rank $k $.$\hfill \blacksquare$\\
\end{prop}

Then the following Corollary is obtained as Corollary 3.3.10 in \cite{[B.24]} since the analog of Lemma 3.3.9 in {\it loc. cit.}  is immediate.

\begin{cor}\label{connection 2}
For each $\lambda \in \C \setminus \mathbb{N}$ there exists a meromorphic integrable connection $\nabla_{\mu, \lambda} : \mathcal{O}_N^k \to \frac{1}{\sigma_k.\Delta}.\mathcal{O}_N^k\otimes \Omega^1_N$ with a simple pole on the reduced hyper-surface $\{\sigma_k\Delta(\sigma) = 0 \} \subset N$ such that the restriction of  $\mathcal{N}_{\mu, \lambda}$ to the Stein (in fact affine)  open set  \ $U := \{\sigma_k\Delta(\sigma) \not= 0 \}$ is isomorphic to the 
$\mathcal{D}_U$-module defined by $(\mathcal{O}_N^k, \nabla_{\mu, \lambda})$. Moreover, this isomorphism is the restriction of an injective $\mathcal{D}_N$-linear map 
$$\mathcal{N}_{\mu, \lambda} \to \big(\mathcal{O}_N^k(*\sigma_k\Delta(\sigma)), \nabla_{\mu, \lambda}\big).$$
This implies that for $\lambda \in \C \setminus \mathbb{N}$ the $D_N$-module $\mathcal{N}_{\mu, \lambda}$ has no torsion.$\hfill\blacksquare$
\end{cor}

\begin{thm}\label{minimal} {\bf [Minimality of $\mathcal{N}_{\mu, \lambda}$]}
Let  $(\mu, \lambda) \in \C^2 $ such that the complex numbers $\lambda, \mu$  are not in $\mathbb{Z}$.
 Then $\mathcal{N}_{\mu, \lambda}$ is the minimal extension of the meromorphic connection given by the rank $k$ vector bundle  $(\mathcal{N}_{\mu, \lambda}(1), \nabla_{\mu, \lambda})$.
\end{thm}

\parag{Remark} This meromorphic connection is exhibited for $\lambda \not\in \mathbb{Z}$ and $\lambda \not\in -k\mu + \mathbb{Z}$ in the last  paragraph of this section.

\parag{Proof} To see that $\mathcal{N}_{\mu, \lambda}$ is the minimal extension of the simple pole meromorphic  connection $\big(\mathcal{N}_{\mu, \lambda}(1), \nabla_{\mu, \lambda}\big)$ it is enough to prove that $\mathcal{N}_{\mu, \lambda}$ has no torsion, and this is given by Corollary \ref{connection 2}, since $\lambda \not\in  \mathbb{Z}$,   and no {\it co-torsion}, that is to say that there is no non trivial coherent left ideal $\mathcal{K}$ in $\mathcal{D}_N$  containing $\mathcal{J}_{\mu,\lambda}$ and generically equal to $\mathcal{J}_{\mu, \lambda}$ on $N$. Such an ideal defines a  holonomic quotient $Q$ of $\mathcal{N}_{\mu, \lambda}$ which is supported in a closed analytic subset $S$ of $N$ with empty interior in $N$. As we may apply Corollary \ref{co-torsion 2}  it is enough to show that near the generic points of $\{\sigma_k\Delta(\sigma) = 0\}$ such an ideal $\mathcal{K}$ is equal to $\mathcal{J}_{\mu, \lambda} $ or to $\mathcal{D}_N$.\\
Near the generic point of $\{\sigma_k = 0\}$ we have $\Delta \not= 0$ and we may use a local isomorphism of $N$ given by a holomorphic section of the quotient map
 $$quot : M = \C^k \to \C^k\big/\mathfrak{S}_k = N.$$
 Via such an isomorphism the ideal  $\mathcal{J}_{\mu, \lambda}$ contains the left ideal $\mathcal{J}^0_{\mu, \lambda}$ in $D_M$ 
 generated by the partial differential operators
 $$\nabla_{i, j} := (z_i - z_j)\frac{\partial^2}{\partial z_i \partial z_j} - \mu(\partial_{z_i} - \partial_{z_j})$$
  for $i \not= j \in [1, k]$ and $\sum_{j=1}^k z_j\frac{\partial}{\partial z_j} - \lambda$, thanks to Lemma \ref{on M} and then  Proposition  \ref{no co-torsion 1}  below allow to conclude this case
 (no $\sigma_k$ co-torsion) using Theorem \ref{isom 2}.\\
The other case (no $\Delta$ co-torsion), that is to say  the vanishing of such a quotient of $\mathcal{N}_{\mu, \lambda}$ near  the generic point of $ \{ \Delta = 0 \}$ is solved in the next paragraph.

\begin{lemma}\label{on M}
The left ideal in $D_M$ which is the annihilator for each $a$  with $\vert a\vert$ large enough compare to $z_1, \dots, z_k$,  of the functions $\big(\prod_{j=1}^k (a - z_j)\big)^{-\mu}$ contains the differential operators
$$ \nabla_{i, j} := (z_i - z_j)\frac{\partial^2}{\partial z_i \partial z_j} - \mu(\partial_{z_i} - \partial_{z_j}), \quad \forall  i, j \quad 1 \leq i < j \leq k .$$
Then the pull-back in $D_M$ of the left ideal of $\mathcal{I}_\mu$  by the (local)  isomorphism given by  a local  holomorphic section of $q$ outside $\{\Delta = 0\}$   contains  these partial differential operators.
\end{lemma}

\parag{Proof} Let us first compute on the open set $\{\delta:= \prod_{1 \leq i < j \leq k} (z_i - z_j) \not= 0\}$ in $M := \C^k$  
$$ \frac{\partial}{\partial z_i}[P_\sigma(a)] = -\frac{P_\sigma(a)}{a - z_i}  $$
and write, for instance when $\vert a\vert$ is large enough compared with $\vert\vert \sigma \vert\vert$:
$$ P_\sigma(a)^{-\mu} = \exp[-\mu Log(P_\sigma(a))] .$$
We obtain 
$$  \frac{\partial}{\partial z_i}[P_\sigma(a)^{-\mu}]  = \mu\frac{1}{a - z_i}P_\sigma(a)^{-\mu} $$
and also
$$ \frac{\partial^2}{\partial z_i\partial z_j}[P_\sigma(a)^{-\mu}] = \mu^2\frac{1}{(a - z_i)(a - z_j)}P_\sigma(a)^{-\mu} .$$
Then the equality
$$ \frac{1}{(a - z_i)} - \frac{1}{(a - z_j)} = \frac{z_i - z_j}{(a - z_i)(a - z_j)}$$
allows to conclude.$\hfill\blacksquare$\\

\begin{prop}\label{no co-torsion 1}
Let $\mathcal{J}^0_{\mu, \lambda}$ for $\lambda + (k-1)\mu \not\in -\mathbb{N}^*$ be the ideal in $D_{\C^k}$ generated by the differential operators  
$$\nabla_{i, j} := (z_i - z_j)\frac{\partial^2}{\partial z_i\partial z_j} - \mu(\partial_{z_i} -\partial_{z_j})$$
for $1 \leq i < j \leq k  $ and
 $$U_0 - \lambda :=  \sum_{h=1}^k z_h\frac{\partial}{\partial z_h} - \lambda.$$
Let assume that $Q$ is a quotient of the $D_{\C^k}-$module $\mathcal{N}^0_{\mu, \lambda} := D_{\C^k}\big/\mathcal{J}^0_{\mu, \lambda}$  in a neighborhood $U$ of the point $(z^0_1, \dots, z^0_k)$ in $\C^k$ where $z_1^0 = 0$ and $z_i \not= z_j$ for $1 \leq i < j \leq k$, with support in $\{z_1 = 0 \}$. Then $Q = 0$.
\end{prop}

\parag{Proof}Assume that $Q \not= 0$ Then $Q = D_U\big/\mathcal{K}$ where $\mathcal{K}$ is a left ideal in $D_U$ such that $\mathcal{J}^0_{\mu, \lambda} \subsetneq \mathcal{K} \subsetneq D$. Then restricting the open neighborhood $U$ of $z^0$ if necessary, there exists a positive integer $n$ such that $z_1^n $ belongs\footnote{The class of $1$ in $Q$ is of $z_1-$torsion !} to $\mathcal{K}$. Then the following lemma allows to conclude the proof when $\lambda+ (k-1)\mu  \not\in -\mathbb{N}^*$.

\begin{lemma}\label{Analogue du Lemme 3.3.12}
Let $K$ be a left  ideal in $D_M$ containing   strictly $\mathcal{J}^0_{\mu, \lambda}$ and  such that the support of $D_M/K$ is in $\{z_1 = 0\}$. Then $K = D_M$ near the generic point of $\{z_1 = 0 \}$ when $\lambda+ (k-1)\mu  \not\in -\mathbb{N}^*$.
\end{lemma}

\parag{Proof} Assume that $z_1^n$ is in $K$ for some $n \in \mathbb{N}^*$. Then using the fact that  $\sum_{h=1}^k z_h\partial_h - \lambda$ is in $K$ we obtain
$$ \partial_1 z_1^n = z_1^n\partial_1 + nz_1^{n-1} =  (n + \lambda) z_1^{n-1} - z_1^{n-1}\sum_{h=2}^k z_h\partial_h \in K .$$
Since $\nabla_{1, j} = \partial_1\partial_j - \frac{\mu}{z_1 - z_j}(\partial_1 - \partial_j) $ we have, since $\nabla_{1, j} $ is in $K$:
$$\nabla_{1, j}z_1^n - z_1^n\nabla_{1, j}  = \frac{-n\mu}{z_1 - z_j}z_1^{n-1} + n z_1^{n-1}\partial_j  \in K $$
Multiplying this on the left by $z_j$  and summing for $j \in [2, k]$ gives
$$ nz_1^{n-1}\sum_{h=2}^k z_j\partial_j - n\mu z_1^{n-1} \sum_{j=2}^k \frac{z_j}{z_1 - z_j}  \in K $$
and using again the fact that  $z_1^{n-1}\sum_{j=2}^k z_j\partial_j = (n+ \lambda)z_1^{n-1}$ modulo  $K$ proved above, we obtain
\begin{align*}
& nz_1^{n-1}(n + \lambda) - n\mu z_1^{n-1} \sum_{j=2}^k (\frac{z_1}{z_1 - z_j} - 1) \in K \\
& n\big(n + \lambda   +  (k-1)\mu\big)z_1^{n-1} - n\mu z_1^n \sum_{j=2}^k (\frac{1}{z_1 - z_j})  \in K 
\end{align*} 
 and using the fact that $z_1^n \in K$ and the non vanishing of  $z_1 - z_j, j \in [2, k]$, at the generic point in $\{z_1 = 0\}$, we conclude 
that $$(\lambda + (k-1)\mu + n) z_1^{n-1} \in K.$$ 
So $z_1^{n-1}$ is in $K$ since $\lambda + (k-1)\mu$ is not in $-\mathbb{N}^*$. So if $Q := D_M/K$ has some $z_1$-torsion near the generic point of the hyper-surface 
$\{z_1 = 0 \}$
we conclude that $Q = 0$ near this generic point.$\hfill \blacksquare$\\

\parag{No $\sigma_k$ co-torsion}  We assume that $\mu$ and $\lambda$ are not in $\mathbb{Z}$. Since there is no $\sigma_k$ co-torsion for $\mathcal{N}_{\mu, \lambda}$ when  $\lambda = -(k-1)\mu$. Then   Theorem \ref{isom 2} implies that for any $\mu$ and $\lambda$  not in $\mathbb{Z}$ the $D_N$-module $\mathcal{N}_{\mu, \lambda}$  has no  $\sigma_k$ co-torsion, concluding this part of  the proof of Theorem \ref{minimal}.\\

 We shall now complete the proof of Theorem \ref{minimal} by proving that a coherent quotient of $\mathcal{M}_\mu$ which is supported by the hyper-surface $\{\Delta = 0\}$ must vanishes.\\
 Recall that for obtaining this result  it is enough, thanks to Corollary \ref{co-torsion 2},  to prove that such a  coherent quotient vanishes around the generic point of this hyper-surface.

 \parag{End of proof of Theorem \ref{minimal}: no $\Delta$ co-torsion} \begin{itemize}
 \item  Near a generic point $\sigma^0$ of the hyper-surface $\{ \Delta(\sigma) = 0 \}$ in $N$ we shall use local coordinates $\tau_1, \tau_2, z_3, \dots, z_k$ where $\tau_1, \tau_2$ are the elementary symmetric functions of the two roots of $P(\sigma)$ which are near the double root  $ z_1 = z_2$ of $P(\sigma^0)$ and where $z_3, \dots, z_k$ are the $k-2$ other roots of $P(\sigma)$ (see Lemma 2.3.3 in {\it loc. cit.}). Then we have a local chart $(\tau, \theta) : U \to U_2\times U_{k-2} $ of an open neighborhood $U$ of $\sigma^0$  in $N$ to a product of polydisc in $\C^2$ and $\C^{k-2}$ which allows to consider $D_{U_2}$ and $D_{U_{k-2}}$ as embedded in $D_U$ as commuting sub-algebras. Now define $\Delta_2 := \tau_1^2 - 4 \tau_2$ and  remark that there is an invertible holomorphic function $I$ on $U$ such that $\Delta(\sigma) =  I(\sigma)\Delta_2(\tau)$ on $U$. Define also on $U$ the differential operator
 $$T^2_\mu := \partial_{\tau_1}^2 + \tau_1\partial_{\tau_1}\partial_{\tau_2} +  \tau_2\partial_{\tau_2}^2 + (\mu + 1) \partial_{\tau_2}   $$

 \begin{lemma}\label{co-torsion delta}
Let $\sigma^0 := q(z_1^0, z_1^0, z_3^0, \dots, z_k^0)\in N$ be a generic point in the hyper-surface  $\{\Delta = 0 \}$,  where the complex numbers $z_1^0, z_3^0, \dots, z_k^0$ are pair-wise distinct. Then we  define a  local coordinate system  on $N$ near $\sigma^0$ given by $\tau_1 := z_1 + z_2, \tau_2 := z_1z_2$, and with  $ z_3, \dots, z_k$ where $z_1, z_2$ are the roots of $P_\sigma$ near to the double root $z_1^0 = z_2^0$  of $P_{\sigma^0}$ and $z_j$ the root of $P_\sigma$ near $z^0_j$ for $j \in [3, k]$. \\
Then the partial differential operator $T^2_\mu := \partial_1^2 + \partial_2(\tau_1\partial_1 + \tau_2\partial_2 + \mu) $ belongs to the ideal $\mathcal{I}_\mu$ where $\partial_i$ is the partial derivative in $\tau_i$ for $i = 1, 2$.
\end{lemma}

\parag{Proof} It is enough to show that $T^2_\mu$ annihilates the holomorphic  function $\big(\prod_{j=1}^k (1 - z_j/a)\big)^{-\mu}$ for any $a $ large enough compared to $\vert \sigma^0\vert$. So it is enough to show that for $a$ large enough  $T^2_\mu[(1 - \tau_1/a + \tau_2/a^2)^{-\mu}] = 0 $.\\
We have
\begin{align*}
& \partial_1(1 - \tau_1/a + \tau_2/a^2) = -1/a,\\
& \partial_2(1 - \tau_1/a + \tau_2/a^2) = 1/a^2 \quad {\rm and \ so}\\
& \partial_1^2[(1 - \tau_1/a + \tau_2/a^2)^{-\mu}] = \frac{\mu(\mu + 1)}{a^2}(1 - \tau_1/a + \tau_2/a^2)^{-\mu - 2} \quad {\rm and} \\
& \big(\tau_1\partial_1 + \tau_2\partial_2 + \mu\big) [(1 - \tau_1/a + \tau_2/a^2)^{-\mu}] = \mu(1 - \tau_1/a + \tau_2/a^2)^{-\mu - 1} \quad {\rm and} \\
& \partial_2 \big(\tau_1\partial_1 + \tau_2\partial_2 + \mu\big) [(1 - \tau_1/a + \tau_2/a^2)^{-\mu}] =  -\frac{\mu(\mu + 1}{a^2}(1 - \tau_1/a + \tau_2/a^2)^{-\mu - 2}
\end{align*}
The conclusion follows.$\hfill \blacksquare$

\parag{Remark}Since we have $P(\sigma, a) = (a^2 - \tau_1a + \tau_2)\prod_{j=3}^k (a - z_j)$ in our local coordinates we may write for each $\mu \in \C$ and each $q \in \mathbb{N}$:
 $$ \Phi^{(k)}_{\mu, q}(\sigma) = \sum_{j+j' = q} \Phi^{(2)}_{\mu, j}(\tau)\Phi^{(k-2)}_{\mu, j'}(z_3, \dots, z_k) $$
 where we define 
 $$\big(\prod_{j=1}^2 (1 - z_j/a)\big)^{-\mu} =  \Big(\frac{a^2}{a^2 - \tau_1a + \tau_2}\Big)^\mu := \sum_{j=0}^\infty \Phi^{(2)}_{\mu, j}(\tau)a^{-j} $$ 
 and
 $$ \big(\prod_{h=3}^k (1 - z_h/a)\big)^{-\mu} = \sum_{j'= 0}^\infty \Phi^{(k-2)}_{\mu, j'}(z_3, \dots, z_k)a^{-j'} \qquad \qquad \square$$

 \item Lemma 2.3.3 in {\it loc. cit.} gives the identities
 \begin{equation*}
  (T^2_\mu - 2n\partial_{\tau_2})\Delta_2^n = \Delta_2^nT^2_\mu + 2n(2n+1-\mu)\Delta_2^{n-1} \quad \forall n \in \mathbb{N}^*  \tag{$@$}
  \end{equation*}
 \item Assume now that  $\mathcal{Q}$ is a non zero   coherent quotient of $\mathcal{M}_\mu$ with support in $\{\Delta = 0 \}$ and that $\mu$ is not in $2\mathbb{N} + 1$. Note $\mathcal{K}$ the left ideal such that $\mathcal{Q} = D_N/\mathcal{K}$. Then $\mathcal{K}$ contains $\mathcal{I}_\mu$ and so $T^2_\mu$. 
 There exists a smallest integer $n \in \mathbb{N}$ such that $\Delta_2^n$ kills the class of $1$ in  $\mathcal{Q}$ near the given generic point $\sigma^0$ in $U$. Then formula $(@)$ implies that $\Delta_2^{n-1}$ also kills this generator near $\sigma^0$. So we see that  $n = 0$ and so $\mathcal{Q} = 0$  since $\mathcal{Q}$ is coherent and with support in $\{\Delta = 0 \}$.$\hfill\blacksquare$\\
\end{itemize}

  \subsection{The cases $ \lambda \not\in - k\mu + \mathbb{Z}, \lambda \not\in \mathbb{Z}$}

The aim of this paragraph is to exhibit the meromorphic connection whose minimal extension is the $D_N$-module $\mathcal{N}_{\mu, \lambda}$ for $ \lambda \not\in - k\mu + \mathbb{Z}, \lambda \not\in \mathbb{Z}$\\
We shall use the following elementary Lemma.

\begin{lemma}\label{5/3/26 1} Let $M = \begin{pmatrix} T & Z \\ X & Y \end{pmatrix}$ be a $(p+q, p+q)$ matrix where $T$ is an invertible $(p, p)$ matrix and where $Y$
is a $(q, q)$-matrix. Then $M$ is invertible if and only if the $(q, q)$-matrix $Y - XT^{-1}Z$ is invertible. \\
\end{lemma}

\parag{Proof} Assume that $M(u, v) = 0$ for some  $u \in K^p$ and some $v \in K^q$. Then this is equivalent to
\begin{align*}
& T(u) + Z(v) = 0 \quad {\rm and} \\
& X(u) + Y(v) = 0 
\end{align*}
and this is equivalent to $u = -T^{-1}Z(v)$ and $Y(v) - XT^{-1}Z(v) = 0$. So the existence of $(u, v)$ in the kernel of $M$ implies that $v$ is in the kernel of
$Y - XT^{-1}R$. Moreover, if $v = 0$ we have $u = 0$. So if $(u, v) \not= (0, 0)$ then $v \not= 0$. Conversly the existence of $v \not= 0$ in the kernel of
$Y - XT^{-1}Z$ implies the existence of $(-T^{-1}(Z(v), v) \not= (0, 0)$ in the kernel of $M$.$\hfill \blacksquare$\\

\begin{cor}\label{5/3/26 2}
Under the assumption of the previous Lemma we have
$$ det(M) = det(T)det(Y - XT^{-1}Z).$$
\end{cor}

\parag{Warning} For $det(T) = 0$ the determinant of $M$ may be non zero ! For instance for  $M = \begin{pmatrix} 0 & 1\\ 1 & 1\end{pmatrix}$ (where $p = q = 1$) satisfies
$det(T) = 0$ and $det(M) = -1$.\\ 
For $det(T)= 0$ the term $det(Y - XT^{-1}Z)$ in our formula is not defined. \\
To  understand better what happens
look at the case $\begin{pmatrix} \varepsilon & 1\\ 1 & 1\end{pmatrix}$ when $\varepsilon \to 0$.

\parag{Proof}
 Consider the rational function  $G := det(Y - XT^{-1}Z)$ on $\mathbb{C}^{p+q}$. It is holomorphic
 on the Zariski open set $\{det(T) \not= 0 \}$ which is homogeneous. The function $G$ is homogeneous of degree $q$ on this open set and vanishes only at point where
 $det(M) = 0$. So the function $det(T)G/det(M)$ is holomorphic on this open set and is homogeneous of degree $0$. On each line through the origin which is  not contained
  in $\{det(T) = 0\}$ it is constant. This implies that it is a constant function since it differential vanishes at each point on this connected dense open set.$\hfill \blacksquare$\\

Now we are ready to describe the meromorphic connection associated to the  $D_N$-module $\mathcal{N}_{\mu, \lambda}$ 
 for $\lambda \not\in \mathbb{Z} , \lambda  \not\in - k\mu + \mathbb{Z}$. 
Recall that this $D_N$-module is the quotient of $D_N$ by the left ideal $\mathcal{I}_\mu + D_N(U_0 - \lambda)$. 

\parag{Notations} Let $N := \mathbb{C}^k$ with coordinates $\sigma_1, \dots, \sigma_k$ and denote by $\partial_1, \dots, \partial_k$ the corresponding derivations on $\mathcal{O}_N$.
We denote by $U_0 := \sum_{h=1}^k h\sigma_h\partial_h$ and  by. $U_{-1} := \sum_{h=0}^{k-1} (k-h)\sigma_h\partial_{h+1}$ with the convention $\sigma_0 \equiv 1$

Working modulo $\mathcal{A}$ the class of $\partial_h\partial_p$ depends only on $h + p$ and will be denoted by  $y_{h+p}$. So in $\mathcal{W} := D_N\big/ \mathcal{A}$ the $\mathcal{O}_N$-sub-module $\mathcal{W}(2)\big/ \mathcal{W}(1)$ is free of basis $y_2, \dots, y_{2k}$. We denote by $B$ and $C$ the sub-modules   generated respectively by $y_2, \dots, y_k$ and $y_{k+1}, \dots, y_{2k}$. Their ranks are respectively $(k-1)$ and $k$.\\
To study the class induced by differential operators of order at most $2$ in $\mathcal{N}_{\mu, \lambda}$ we may consider only the $\mathcal{O}_N$ modules generated by 
$A := \sum_{h=1}^k \mathcal{O}_N\partial_h, B$ and $C$ thanks to the relations (modulo our ideal $ \mathcal{I}_\mu + D_N(U_0 - \lambda) $)
$$ 0 \not=  \lambda = U_0 = \sum_{h=1}^k h\sigma_h\partial_h .$$
Then the generators $T^m_\mu$ in $\mathcal{I}_\mu$  give the $(k-1)$ relations
$$ y_m + \sum_{h=1}^k \sigma_hy_{m+h} = -(\mu + 1)\partial_m \quad m \in [2, k]. $$
They will be denoted by $TB + ZC = -(\mu + 1)A'$ where $A'$ is the $(k-1)$-column $^t(\partial_2, \dots, \partial_k)$, $B $ is the $(k-1)$-column $^t(y_2, \dots, y_k)$, $C$ is the $k$-column $^t(y_{k+1}, \dots, y_{2k})$ and where $T$ and $Z$ are the following
 $(k-1, k-1)$ and $(k-1, k)$ matrices
$$ T := \begin{pmatrix} 1 & \sigma_1 & \dots & \dots & \sigma_{k-2} \\ 0 & 1 & \sigma_1& \dots & \sigma_{k-3}\\ \dots & \dots & 1 &\dots & \dots \\ \dots & \dots & \dots &\dots & \sigma_1 \\ 0 & \dots & \dots & 0 & 1 \end{pmatrix} $$
$$ Z = \begin{pmatrix} \sigma_{k-1} & \sigma_k & 0 & \dots & \dots & 0 \\ \sigma_{k-2} & \sigma_{k-1} & \sigma_k & 0 &  \dots & 0\\ 
\dots & \dots & \dots & \dots & \dots & \dots \\ \dots & \dots & \dots & \sigma_k & \dots & 0 \\ \sigma_1 & \dots & \dots & \dots &  \sigma_{k-1} &\sigma_k \end{pmatrix} .$$

Now the relations given by $\partial_h(U_0 - \lambda)$ will be denoted by $XB + YC = UA$ where $X, Y, U$ are matrices of respective sizes $(k, k-1), (k, k)$ and $(k, k)$ and where 
$A$ is the $k$-column $^t(\partial_1, \dots, \partial_k)$. The matrices $X, Y, U$  are defined as follows:
$$ X = \begin{pmatrix}  \sigma_1 & 2\sigma_2 & \dots &\dots & (k-1)\sigma_{k-1} \\ 0 &  \sigma_1 & 2\sigma_2 & \dots & \dots \\ \dots &  \dots & \dots & \dots & \dots \\ \dots & \dots & \dots & \dots & 2\sigma_2 \\  \dots & \dots & \dots & \dots & \sigma_1  \\0 &  \dots & \dots & \dots & 0  \end{pmatrix} $$
$$ Y = \begin{pmatrix}  k\sigma_k & 0 & \dots & \dots & \dots & 0 \\ (k-1)\sigma_{k-1}  & k\sigma_k &  0 & \dots &\dots &  0 \\\dots & \dots & \dots & \dots & \dots & \dots \\  \dots & \dots & \dots & \dots & \dots & \dots \\ \dots & \dots & \dots & \dots & k\sigma_k & 0 \\ \sigma_1 & \dots & \dots & \dots & (k-1)\sigma_{k-1} &k \sigma_k \end{pmatrix} $$

and $U$ is the diagonal matrix with diagonal elements  $(\lambda - p)$ for $p \in [1, k]$.\\
We may write all the relations as
$$ \begin{pmatrix} T & Z \\ X & Y\end{pmatrix}\begin{pmatrix}  B \\ C \end{pmatrix} = \begin{pmatrix} -(\mu+1)A' \\ UA\end{pmatrix} $$
The first relation   $TB + ZC = -(\mu + 1)A'$  allows to compute 
$$B = -T^{-1}ZC -(\mu+1)T^{-1}A'$$
 and so the second relations compute $C$ in term of $A$:
$$ (Y - XT^{-1}Z)C = (\mu+1)XT^{-1}A' + UA $$
where, thanks to the facts that $ det \begin{vmatrix}T & Z \\ X & Y\end{vmatrix} = \sigma_k\Delta$ (see \cite{[B.24]} Lemma 3.3.8 ), the equality  $\det(T) = 1$  and   to Corollary \ref{5/3/26 2}, we know that
$$ det(Y - XT^{-1}Z) = \sigma_k\Delta $$
and we obtain that   the $D_N$-module structure on  $\mathcal{N}_{\mu, \lambda}$ define on $\mathcal{N}_{\mu, \lambda}(1)$  which is a frre $\mathcal{O}_N$-module with basis $\partial_1, \dots, \partial_k$,  a connection with a simple pole on $\{\sigma_k\Delta(\sigma) = 0 \}$

\parag{Example} The computation of $\mathcal{N}_{\mu, \lambda}$ for $k = 2$ and $\mu, \lambda$  are not in $\mathbb{Z}$.\\

The $\mathcal{O}_N$-generators of $\mathcal{N}_{\mu, \lambda}(2)$ are $\partial_1, \partial_2, y_2, y_3, y_4$ where $y_p$ is the class of $\partial_i\partial_j$ with $i+j = p$ since $\lambda \not= 0$.\\
Then we have the relations, where we define $e_1 := \partial_1$ and $e_2 := \partial_2$:
\begin{align*}
& T^2_\mu = 0 \quad {\rm which \ gives} \quad y_2 + \sigma_1y_3 + \sigma_2y_4 = -(\mu+1)e_2\\
& \partial_1(U_0 - \lambda) = 0 \quad {\rm which \ gives} \quad  \sigma_1y_2 + 2\sigma_2y_3 = (\lambda - 1)e_1 \\
& \partial_2(U_0 - \lambda) = 0 \quad {\rm which \ gives} \quad  \sigma_1y_3 + 2\sigma_2y_4 = (\lambda - 2)e_2 .
\end{align*}

The first two relations implies
$$ (\sigma_1^2 - 2\sigma_2)y_3 + \sigma_1\sigma_2y_4 =  -(\mu + 1)\sigma_1e_2 - (\lambda - 1)e_1 .$$
Combined with the third relation we obtain
$$ -\Delta y_3 = 2(\lambda - 1)e_1 + (2\mu + \lambda)\sigma_1e_2 $$
and 
$$ \sigma_2\Delta y_4 = (\lambda - 1)\sigma_1e_1 + \big[(\mu + \lambda - 1)\sigma_1^2 - 2(\lambda - 2)\sigma_2\big]e_2 $$
which may be written
$$ \sigma_2\Delta y_4 = (\lambda - 1)\sigma_1e_1 + \big[(\mu + \lambda - 1)\Delta  + 2(2\mu + \lambda)\sigma_2\big]e_2       $$
and finally
$$ \Delta y_2 = (\lambda - 1)\sigma_1e_1 + 2(2\mu + \lambda)\sigma_2e_2.$$
Since $\partial_1e_1 = y_2, \partial_2e_1 = \partial_1e_2 = y_3$ and $\partial_2 e_2 = y_4$ the formulas above determines the meromorphic connection of the rank $2$
trivial vector bundle with basis $e_1, e_2$ whose minimal extension is the $D_N$-module $\mathcal{N}_{\mu, \lambda}$.\\

\begin{align*}
& \Delta \partial_1e_1 = (\lambda - 1)\sigma_1e_1 + 2(2\mu + \lambda)\sigma_2e_2. \\
& \Delta\partial_2e_1 = - 2(\lambda - 1)e_1 - (2\mu + \lambda)\sigma_1e_2 . \\
& \Delta\partial_1e_2 =  - 2(\lambda - 1)e_1 - (2\mu + \lambda)\sigma_1e_2 . \\
& \sigma_2\Delta \partial_2e_2 = (\lambda - 1)\sigma_1e_1 + \big[(\mu + \lambda - 1)\sigma_1^2 - 2(\lambda - 2)\sigma_2\big]e_2. 
 \end{align*}
 
It is easy to verify directely on these formulas  that we have $(U_0 - \lambda)(e_1) = -e_1$ and $(U_0 - \lambda)(e_2) = -2e_2$ which corresponds the  weight $\lambda$ for the $D_N$-generator  given by  class of $1$  in  this $D_N$-module.

\section{Periodicity in $\mu$}

\subsection{The action of $\square(E+\mu)$ on $\mathcal{M}_{\mu + 1}$}

 First we examine the action of the right product by $(E + \mu)$ on $\mathcal{M}_{\mu+1}$. We shall assume $\mu \not= 0$ in the sequel.

   \begin{prop}\label{period mu}
   For each $\mu \in \C^*$  we have an injective left $D_N$-linear map 
    $$\square (E + \mu) : \mathcal{M}_{\mu+ 1} \to \mathcal{M}_\mu$$
    given by right product by $E + \mu$.
   \end{prop}
   
   \parag{Proof} Thanks to point 4 in Proposition \ref{31/05 suite 1}   the right product by $E + \mu$ sends $\mathcal{M}_{\mu + 1}$ to $\mathcal{M}_\mu$. But the fact that the symbol of $E$ does not vanish identically on  any non empty subset in $Z$ implies more: assume that $\Pi(E+ \mu) $ is in $\mathcal{I}_{\mu, \sigma^0}$ and that $\pi$ has minimal order inside the subset of  germs at $\sigma^0$  of partial differential operators   in $D_N \setminus \mathcal{I}_{\mu + 1}$ with this property (assuming that this subset  is non empty). Then the product of symbols  $s(\Pi)s(E)$ must vanish on $Z$ and so $s(\Pi)$ vanishes on $Z$ and there exists $P \in \mathcal{I}_{\mu+1, \sigma^0}$ with $s(P) = s(\Pi)$. Now $\Pi - P$  is zero or has an order strictly less than the order of $\Pi$. But $(\Pi - P)(E + \mu)$ is in $\mathcal{I}_{\mu, \sigma^0}$ and so $\Pi =  P$. \\
  We conclude that the left ideal in $D_N$ of partial differential operators $\Pi$ such that $\Pi(E + \mu)$ belongs to $\mathcal{I}_\mu$ is equal to $\mathcal{I}_{\mu + 1}$. So the map $\square (E + \mu)$ is injective.\\
  
  Note that the characteristic variety of $\mathcal{M}_\mu\big/ \mathcal{M}_{\mu+1}(E+\mu)$ is the co-normal to the hyper-surface $\{\sigma_k = 0\}$ (with multiplicity $k$, see \cite{[B.24]} Section 3.2). The following Proposition explicits the localization by $\sigma_k^{-1}$ of the co-image of the map $\square(E + \mu)$ and its Corollary shows that this co-image is in fact already localized in $\sigma_k^{-1}$.
  
  \parag{Remark} The image of the map $\square (E + \mu)$ in $\mathcal{M}_\mu$ is isomorphic to the quotient of $D_N$ by the left ideal  $L_\mu$ generated by
   \begin{equation*}
     \mathcal{A}, \  \sum_{m=1}^{k-1} D_N\partial_1\partial_m  \quad {\rm and} \quad  E + \mu  \tag{$L_\mu$}
     \end{equation*}
   
   since $T^{m+1}_\mu = \partial_1\partial_m + \partial_{m+1}(E + \mu)$ for each $m \in [1, k-1]$.\\
   
   \begin{prop}\label{14/12/25}
   The localization in $\sigma_k$ of the  quotient $D_N\big/L_\mu$ is isomorphic  to $(\mathcal{O}_N[\sigma_k^{-1}])^k$ as a $\mathcal{O}_N[\sigma_k^{-1}]$-module with basis $\partial_1, \dots, \partial_k$ where the action of $D_N$ is
   given by the equalities
  $$\partial_k \begin{pmatrix} \partial_1 \\ \dots \\ \partial_k \end{pmatrix} = - (\mu+1) M^{-1}\begin{pmatrix} \partial_1 \\ \dots \\ \partial_k \end{pmatrix} \qquad\qquad  \qquad\qquad  \qquad\qquad (@)$$
  and 
  \begin{align*}
  & \partial_p\partial_q = \partial_k\partial_h \quad {\rm for} \quad p + q = k+h, \  \forall h \in [1, k]       \tag{@@} \\
  & \partial_p\partial_q = 0 \quad {\rm for} \quad p + q \leq k.\tag{@@@}
  \end{align*}
  where 
  $$ M = \begin{pmatrix} \sigma_k & 0 & \dots & \dots & 0 \\ \sigma_{k-1} & \sigma_k & 0 & \dots & 0 \\ \dots &\dots& \dots &\dots &\dots \\ \sigma_2 & \dots & \dots & \sigma_k &0 \\
  \sigma_1 & \dots & \dots & \sigma_{k-1} & \sigma_k \end{pmatrix} $$
   \end{prop}
   \bigskip

   \parag{Proof} Since  $L_\mu$ contains $\mathcal{A}$ and $\partial_1\partial_m$ for each $m \in [1, k-1]$, $\partial_p\partial_q$ is in $L_\mu$ for $p+q \in [1, k]$. So we have for each 
   $h \in [1, k]$
   \begin{equation*}
     \sum_{p+h \in [k+1, 2k]}  \sigma_p\partial_p\partial_h + (\mu+1)\partial_h = \partial_h(E + \mu)  = 0 \quad {\rm (modulo)} \ L_\mu .\tag{M}
     \end{equation*}
   The matrix $M$ defined in our statement is invertible in $\mathcal{O}_N[\sigma_k^{-1}]$ and so $ (\partial_1, \dots, \partial_k)$ generates $\mathcal{O}_N[\sigma_k^{-1}]\otimes_{\mathcal{O}_N} \big(\mathcal{M}_\mu\big/\mathcal{M}_{\mu+1}(E+\mu)\big)$ for $\mu \not= -1$ as a $\mathcal{O}_N[\sigma_k^{-1}]$-module. It is a basis of this module because $\mathcal{I}_\mu$ (and then $L_\mu$\footnote{Since the symbol of $E$ does not vanish on any non empty open subset of $Z$.})  contains no partial differential operator  of order at most $1$  and so $E + \mu = 0$ gives the only relation   between $\partial_1, \dots, \partial_k$ and $1$ (recall that $\mu \not= 0$). So $\partial_1, \dots, \partial_k$ is free on $\mathcal{O}_N$ \\
   The relations $(M)$  in the quotient $\mathcal{M}_\mu\big/\mathcal{M}_{\mu+1}(E+\mu)$  allows to describe the meromorphic connection on this free rank $k$ $\mathcal{O}_N$-module generated by $\partial_1, \dots, \partial_k$ by formulas $(@), (@@)$ and $(@@@)$.$\hfill \blacksquare$\\
    
   \begin{lemma}\label{11/1/26 a} 
   In the left $D_N$-module $\mathcal{M}_\mu\big/\mathcal{M}_{\mu+1}(E+\mu)$ we have the relations
   \begin{equation*}
    M^q\partial_k^q\begin{pmatrix} \partial_1 \\ \dots \\ \partial_k\end{pmatrix} = (-1)^q(\mu+1)\dots(\mu+q)\begin{pmatrix} \partial_1 \\ \dots \\ \partial_k\end{pmatrix} . \tag{$@_q$}
    \end{equation*}
   \end{lemma}
   
   \parag{Proof} For $q = 1$ the relation is given by formula $(@ )$.\\
   So let assume the formula proved for $q \geq 1$ and consider the case $q + 1$. Since we have
   $$ M\partial_kM^q\partial_k^q = M^{q+1}\partial_k^{q+1} + qM^q\partial_k^q$$
   because $M = \sigma_k Id_k + N_k$ with $[\partial_k, N_k] = 0$ implies $[\partial_k, M] = Id_k$ (and then  $[\partial_k, M^q]  = qM^{q-1}$ for each $q$). So  we obtain:
   \begin{align*}
   & M^{q+1}\partial_k^q\begin{pmatrix} \partial_1 \\ \dots \\ \partial_k\end{pmatrix}  = -qM^q\partial_k^q\begin{pmatrix} \partial_1 \\ \dots \\ \partial_k\end{pmatrix} + (-1)^q(\mu+1)\dots(\mu+q)M\partial_k\begin{pmatrix} \partial_1 \\ \dots \\ \partial_k\end{pmatrix} \\
   & \qquad = (-1)^{q+1}(\mu+1)\dots(\mu+q)(\mu+q+1)\begin{pmatrix} \partial_1 \\ \dots \\ \partial_k\end{pmatrix} \\
   \end{align*}
  concluding the proof $\hfill \blacksquare$\\
  
  \begin{cor}\label{11/1/26 b} Assume that $\mu$ is not in $-\mathbb{N}$.  Then for each positive integer $n$ and for each $h \in [1, k]$,  $\sigma_k^{-n}\partial_h$ belongs to $\mathcal{M}_\mu\big/\mathcal{M}_{\mu+1}(E + \mu)$.\\
  So this $D_N$-module is the free $\mathcal{O}_N[\sigma_k^{-1}]$-module with basis $\partial_1,\dots,  \partial_k$ and with the action of $D_N$ described by the formulas $(@), (@@), (@@@)$.
  \end{cor}
  
  As already announced,  this Corollary shows that Proposition \ref{14/12/25} is valid without localization in $\sigma_k$.
  
  \parag{proof} First remark that since $M$ is lower triangular we have 
  $$\sigma_k^q\partial_k^q\partial_1 = (-1)^q(\mu+1)\dots(\mu+q)\partial_1$$
   so that for $\mu \not\in -\mathbb{N}$ we
  have $\sigma_k^{-q}\partial_1$ belongs to $\mathcal{M}_\mu\big/\mathcal{M}_{\mu+1}(E+\mu)$ (in the sense that there exists an element $\xi$  such that $\sigma_k^q \xi = \partial_1$).
  After quotient by the $\mathcal{O}_N[\sigma_k^{-1}]$-sub-module generated by $\partial_1$ (or, if you prefer by the sub-$\mathcal{O}_N$-module generated by the $\partial_k^q\partial_1$ for $q \in \mathbb{N}$), we are in the same situation but replacing the matrix $M$ by  the triangular $(k-1, k-1)$-matrix  obtained by deleting the first row and the first line in $M$. It has still $\sigma_k$ on the diagonal, so still invertible when $\sigma_k \not= 0$. Continuing in this manner we conclude that $\mathcal{M}_\mu\big/\mathcal{M}_{\mu+1}(E+\mu)$ is the $\mathcal{O}_N[\sigma_k^{-1}]$-module generated by $\partial_1, \dots, \partial_k$.\\
  Then $\mathcal{M}_\mu\big/\mathcal{M}_{\mu+1}(E+\mu)$  has no $\sigma_k$-torsion and so  $\partial_1, \dots, \partial_k$  is a $\mathcal{O}_N[\sigma_k^{-1}]$ basis.$\hfill \blacksquare$\\
  
  \parag{Remark} Denote by $[M^{-1}]_{k-p+1}$ the $(k, k)$ matrix  obtained from $M^{-1}$ by replacing the lines $1$ to $k-p$ by zeros. Then the connection of the free, rank $k$, 
   $\mathcal{O}_N$-module generated by $\partial_1, \dots, \partial_k$ is explicitly given by the formulas
  $$ \partial_p\begin{pmatrix}\partial_1 \\ \dots \\ \partial_k \end{pmatrix} = -(\mu+1)[M^{-1}]_{k-p+1}\begin{pmatrix}\partial_1 \\ \dots \\ \partial_k \end{pmatrix} $$
  for each $p \in [1, k]$. Of course, this connection is integrable and has a  pole of order at most equal to $k$   along the hyper-surface $\{\sigma_k = 0 \}$ in $N$.\\
  We conclude that, for $\mu \not\in -\mathbb{N}$,   the $D_N$-module 
  $$\mathcal{M}_\mu \big/\mathcal{M}_\mu (E + \mu) \simeq D_N\big/(\mathcal{I}_\mu + D_N(E + \mu))$$
   is the free, rank $k$,  $\mathcal{O}_N[\sigma_k^{-1}]$-module with basis $\partial_1, \dots, \partial_k$ endowed with the meromorphic connection described above.
  \\
  
  \subsection{The action of $\square(E+\mu)$ on $\mathcal{N}_{\mu+1, \lambda}$}
  
  The goal of this sub-section is the proof of the following theorem

\begin{thm}\label{11/1/26} {\bf[Periodicity in $\mu$]}  Assume that $\mu \not= 0$,  $\lambda$ is not in  $ \mathbb{N}$ and that $\lambda + k\mu \not= 0, -1, \dots, -(k-1)$.
 Then the left $D_N$-linear map $\square (E + \mu): \mathcal{N}_{\mu+1, \lambda} \to \mathcal{N}_{\mu, \lambda}$ is an isomorphism.
 \end{thm}
 
 \parag{Proof} This result  is an easy  consequence of Proposition \ref{period mu} below combined with  Corollary \ref{connection 2} since Proposition \ref{period mu} gives that the map 
  $\square (E + \mu)$ is  surjective. As  the map $\square(E + \mu)$ is an $\mathcal{O}_N$-linear  surjective map between rank $k$ free $\mathcal{O}_U$-modules on the Zariski dense open set $U := \{\sigma_k\Delta(\sigma) \not= 0 \}$ in $N$, its kernel  is supported by an hyper-surface;  so it is a torsion submodule. But Corollary \ref{connection 2} implies that $\mathcal{N}_{\mu + 1, \lambda}$ has no $\mathcal{O}_N$-torsion for $ \lambda$ not in $\mathbb{N}$. The conclusion follows. $\hfill \blacksquare$\\
  
  The cases for which the map $\square(E + \mu)$ is not an isomorphism are examined in Section 5.4.\\
  
  Recall that the ideal $L_\mu$ of $D_N$ is generated by $\mathcal{A}$, $\partial_1\partial_{m-1}$ for $m \in [2, k]$ and $E + \mu$.
 
 \begin{lemma}\label{4/9}
 For $\mu \not= 0$ and $\lambda + k(\mu+ 1) - h \not= 0$ for $h \in [1, k]$  we have    $$ L_\mu + D_N(U_0 - \lambda) = D_N.$$
 \end{lemma}

 \parag{Proof} We first prove that $\partial_h$ belongs to $L_\mu + D_N(U_0 - \lambda)$ for any $h \in [1, k]$.\\
 \begin{align*}
 & {\rm  Since } \qquad \partial_h(E + \mu) = (\mu + 1)\partial_h + \sum_{m=1}^{k-1} \sigma_m\partial_h\partial_m + \sigma_k\partial_h\partial_k \\
 & {\rm and} \qquad  \partial_h(U_0 - \lambda) = (h - \lambda)\partial_h + \sum_{m=1}^{k-1} m\sigma_m\partial_h\partial_m + k\sigma_k\partial_h\partial_k
 \end{align*}
  and $L_\mu$ contains $\partial_1\partial_m $ for each $m \in [1, k-1]$ 
 we obtain that
 $(k\mu + \lambda + k - 1)\partial_1$ is in $L_\mu + D_N(U_0 - \lambda)$. Then for $h = 2$ we obtain in an analogous way
 $$ (k(\mu+1) + \lambda - 2)\partial_2  \in L _\mu + D_N(U_0 - \lambda)$$
 since $\partial_m\partial_2 = \partial_{m+1}\partial_1$ in $L_\mu$ for $m \leq k-1$ and continuing in this way we obtain that $\partial_1, \dots, \partial_k$ are in 
 $L_\mu + D_N(U_0 - \lambda)$. Since $E + \mu$ is in $L_\mu$ with $\mu \not= 0$ we conclude that $L_\mu + D_N(U_0 - \lambda) =  D_N$ for $\mu\not= 0$ and $\lambda \not= h - k(\mu+1)$ for $ h \in [1, k]$.$\hfill\blacksquare$\\
 
 \parag{Remark} For $\mu \not= 0$ and  $h \in [2, k]$ the ideal $L_\mu + D_N(U_0 +k(\mu+1) -  h)$ contains $\partial_1, \dots, \partial_{h-1}$. So the rank of the vector bundle on $N \setminus \{\sigma_k\Delta(\sigma) = 0 \}$ defined by
  $$\mathcal{N}_{\mu, k(\mu + 1) -h}\big/(Im(\square(E + \mu)) = D_N\big/(L_\mu + D_N(U_0 + k\mu - h))$$ 
   is at most $k -h + 1$ and the corresponding kernel of $\square(E+ \mu)$ in $\mathcal{N}_{\mu+1, -k(\mu + 1) - h}$ is at most  equal to $k- h + 1$ and its image has at least rank $h - 1$.$\hfill \square$\\
 
The next proposition completes the proof of Theorem \ref{11/1/26}.
 
 \begin{prop}\label{period mu} 
 For $\mu \not= 0$ and $\lambda + k\mu  \not= 0, -1, \dots, -k+1$ the left  $D_N$-linear map $\square(E + \mu) : \mathcal{N}_{\mu + 1, \lambda} \to \mathcal{N}_{\mu, \lambda}$ is a surjective map. \end{prop}
 
 \parag{Proof} Since  we have $[U_0, E] = 0$, the map $\square(E + \mu) : \mathcal{M}_{\mu + 1} \to \mathcal{M}_\mu$ induces a left $D_N$-linear map
 $\square(E + \mu) : \mathcal{N}_{\mu + 1, \lambda} \to \mathcal{N}_{\mu, \lambda}$. \\
 To prove that this map is surjective remark that the  co-kernel of our map in $\mathcal{N}_{\mu, \lambda}$ is the quotient of $D_N$ by the ideal
$\mathcal{I}_\mu + D_N(U_0 - \lambda) + D_N(E + \mu) $ which is equal to  $L_\mu + D_N(U_0 - \lambda) = D_N$ thanks to Lemma \ref{4/9}.  The conclusion follows $\hfill \blacksquare$\\

 \section{Special cases}
 
 \subsection{Some tools}
 
 We examine in this section some cases in which the previous results do not apply. In this first  paragraph we develop some tools  which will be used later on.\\
 
  The next lemma  shows that, for $\lambda \not\in \mathbb{N}^*$,  the class induces by  $U_{-1}$ in the $D_N$-module  $\mathcal{N}_{\mu, \lambda}$ is not  a  $\mathcal{O}_N$-torsion class.

 \begin{lemma}\label{2/3}
 For any choices of $\mu$ and $\lambda \not\in \mathbb{N}^*$, for any $\sigma \in N$ and any non zero $f \in \mathcal{O}_{N, \sigma}$ we have $fU_{-1} \not\in \mathcal{J}_{\mu, \lambda}$.
 \end{lemma}
 
 \parag{Proof} Assume the lemma is wrong and  assume that for some $\sigma \in N$ and for some non zero $f \in \mathcal{O}_{N, \sigma}$  we have
 $$ fU_{-1} = P + Q(U_0 - \lambda) $$
 where $P$ is in $\mathcal{I}_{\mu, \sigma}$ and $Q \in D_{N, \sigma}$. Assume, moreover that $P$ has the minimal order $\geq 2$ to satisfy such an equality for some given non zero  $f$.
 Since the symbols satisfy $s(P) + s(Q)s(U_0) = 0$ and since $s(P)$ must vanish on $Z$, the fact that $s(U_0)$ does not vanish in any non empty open set in $Z$ forces $s(Q)$ to vanish on $Z$.
 Then there exists $Q_1 \in \mathcal{I}_{\mu, \sigma}$ which satisfies $s(Q_1) = s(Q)$. Now define $P_1 := P -Q_1(U_0 - \lambda)$. Since $\mathcal{I}_\mu(U_0 - \lambda) \subset \mathcal{I}_\mu$ for $\lambda \not\in \mathbb{Z}$\footnote{ We have $(U_0 - \lambda)\Phi_{\mu, q} = (q - \lambda)\Phi_{\mu, q}$ and we use Proposition \ref{31/05 suite 1}.} we have $P_1 \in \mathcal{I}_{\mu, \sigma}$, the order of $P_1$ is strictly less than the order of $P$ and $fU_{-1} = P_1 + (Q - Q_1)(U_0 - \lambda)$. So we obtain the existence of $P \in \mathcal{I}_\mu$ of order at most  $1$ and $Q$ such that $ fU_{-1} = P + Q(U_0 - \lambda)$. But then $P = 0$ (see Lemma  \ref{remplace}) and then $Q$ is in $\mathcal{O}_{N, \sigma}$ and satisfies $fU_{-1} = Q(U_0 - \lambda)$. This is not possible for $\lambda \not= 0$, and the proof is complete in this case.\\
 For $\lambda = 0$ the proof above is still valid because $U_0$ kills $\Phi_{\mu, 0} = 1$ so the reduction to the case $P = 0$ still works; and at the end of the proof we find that
 $fU_{-1} = QU_0$ where $Q$ is in $\mathcal{O}_{N, \sigma}$. This implies the equalities in an open neighborhood of $\sigma$:
 $$ (k-h)f\sigma_h = (h+1)Q\sigma_{h+1}\quad {\rm for } \quad h \in [0, k-1] $$
 which is not possible on a non empty open set in $N$. So the Lemma is also true for $\lambda = 0$.$\hfill \blacksquare$

 \begin{prop}\label{magic 1}
For $h \in [2, k]$ we have the equality
\begin{equation*}
\partial_h(U_0 +k\mu - 1) + \partial_{h-1}U_{-1} = k\mathcal{T}^h_\mu +  \sum_{q=1}^{k-1} (k-q)\sigma_qA_{h-1, q+1} \tag{$E_h(\mu)$}
\end{equation*}
and for $h = 1$ the equality
\begin{equation*}
-\partial_1(U_0+ k\mu - 1) +  (E +\mu)U_{-1} =  \sum_{q=1}^{k-1} (k-q)\sigma_{q}\mathcal{T}^{q+1}_\mu  . \tag{$E_1(\mu)$}
\end{equation*}
\end{prop}

The proof is an obvious consequence of Proposition. 3.1.5   in \cite{[B.24]} which establishes these formulas for $\mu = 0$. The reader may also verify easily  these formulas by a direct computation$\hfill \blacksquare$\\

\begin{prop}\label{magic 2}
For $j \in [2, k]$ we have
\begin{equation*}
\partial_jU_1 +  (\lambda + 1)\partial_{j - 1}  =  - \sum_{h=1}^{k-1} (h+1)\sigma_{h+1}A_{h, j} +  \sigma_1T^j - \partial_{j-1}(U_0 - \lambda) \tag{$F_j$}
\end{equation*}
and for $j = 1$
\begin{equation*}
\partial_1U_1 -  (\lambda + 1)E =  E(U_0 - \lambda) - \sum_{h=1}^{k-1} (h+1)\sigma_{h+1}T^{h+1} \tag{$F_1$}
\end{equation*}
\end{prop}

\parag{Proof} Contrary to the case of the previous formulas, these are not given completely in \cite{[B.24]} (see Lemma 4.2.3), so we detail here the computation.\\
 We have $U_1 = \sigma_1E  - \sum_{h=1}^{k-1} (h+1)\sigma_{h+1}\partial_h $ \  and so, for $j \geq 2$
\begin{align*}
& \partial_jU_1 = \partial_j\sigma_1 E - \sum_{h=1}^{k-1} (h+1)\sigma_{h+1}\partial_j\partial_h  \ - j\partial_{j-1} \\
& \qquad = \sigma_1(T^j - \partial_1\partial_{j-1}) - \sum_{h=1}^{k-1} (h+1)\sigma_{h+1}A_{h, j} - \sum_{h=1}^{k-1} (h+1)\sigma_{h+1}\partial_{h+1}\partial_{j-1} - j\partial_{j-1} \\
& \qquad = \sigma_1T^j - \sigma_1\partial_1\partial_{j-1} - \sum_{h=1}^{k-1} (h+1)\sigma_{h+1}A_{h, j}  - \partial_{j-1}(U_0 - \sigma_1\partial_1)) + (j-1)\partial_{j-1}  - j\partial_{j-1}\\
& \qquad =  -(\lambda + 1)\partial_{j-1}  - \sum_{h=1}^{k-1} (h+1)\sigma_{h+1}A_{h, j} - \partial_{j-1}(U_0 - \lambda) + \sigma_1T^j
\end{align*}
For $j = 1$ we have
\begin{align*}
& \partial_1U_1 = \sigma_1\partial_1E + E - \sum_{h=1}^{k-1} (h+1)\sigma_{h+1}(T^{h+1} - \partial_{h+1}E) \\
& \qquad =  \sigma_1\partial_1E + E + (U_0 - \sigma_1\partial_1)E - \sum_{h=1}^{k-1} (h+1)\sigma_{h+1}T^{h+1} \\
& \qquad = (\lambda + 1)E + E(U_0 - \lambda)  - \sum_{h=1}^{k-1} (h+1)\sigma_{h+1}T^{h+1} 
\end{align*}
since $E$ and $U_0$ commute.$\hfill \blacksquare$\\

\begin{cor}\label{mu magic 2}
For each $\mu \in \mathbb{C}$ we have for $j \in [2, k]$:
\begin{equation*}
 \partial_j(U_1 + \mu\sigma_1)  +  (\lambda + 1)\partial_{j - 1}. =   - \sum_{h=1}^{k-1} (h+1)\sigma_{h+1}A_{h, j} +  \sigma_1T_\mu^j - \partial_{j-1}(U_0 - \lambda) \tag{$F_j(\mu)$}
 \end{equation*}
 and for $j = 1$:
 \begin{equation*}  
 \partial_1(U_1 + \mu\sigma_1) -  (\lambda + 1)(E + \mu)  =  (E + \mu)(U_0 - \lambda)  - \sum_{h=1}^{k-1} (h+1)\sigma_{h+1}T_\mu^{h+1}  \tag{$F_1(\mu)$}
\end{equation*}
\end{cor}

\parag{Proof} For $j \geq 2$ have only to compare the terms in $\mu$ which appear on both sides since the formula is true for $\mu = 0$. The find $\mu\sigma_1\partial_j$ for the left-side and $\sigma_1\mu \partial_j$ for the right-side. This case follows.\\
 For $j = 1$ also  we  have only to compare the terms in $\mu$ which appear on both sides since the formula is true for $\mu = 0$. We find $\mu\partial_1\sigma_1 = \mu + \mu\sigma_1\partial_1$ on the left-side and
$$  \mu(\lambda + 1) + \mu(U_0 - \lambda) - \mu(U_0 - \sigma_1\partial_1)  = \mu\sigma_1\partial_1 + \mu$$
the term  $\mu(U_0 - \sigma_1\partial_1) $ coming from the fact that $T^{h+1}_\mu - T^{h+1} = \mu\partial_{h+1}$ and the definition of $U_0 := \sum_{p=1}^k p\sigma_p\partial_p$.$\hfill \blacksquare$\\

   \subsection{The cases $\lambda = 1 - k\mu$ and $\lambda = -k\mu$   with $ k\mu \not\in \mathbb{Z}$.}

   \begin{prop}\label{torsion 1.1} Assume that $\mu \not\in \mathbb{Z}$ and $k\mu \not= 1$. Then there exists an injective $\mathcal{D}_N$-linear map $\chi : K_\mu \to \mathcal{N}_{\mu, 1- k\mu}$ which sends the class of  $1$ in $K_\mu$ to the class of $U_{-1}$ in $\mathcal{N}_{\mu, 1- k\mu}$. \\
The quotient  $\mathcal{N}^\square_{\mu, 1- k\mu} := \mathcal{N}_{\mu, 1- k\mu}\big/K_\mu$ is  isomorphic, via the map 
$$\square(U_1 + \mu\sigma_1) :  \mathcal{N}_{\mu, 1- k\mu} \to \mathcal{N}_{\mu, - k\mu}$$
 to the  sub-module $\mathcal{N}^\sharp_{\mu, - k\mu}$  of $\mathcal{N}_{\mu, - k\mu}$ which  is  the image of $H_\mu$ in $\mathcal{N}_{\mu, - k\mu}$.\\
 Moreover  $\mathcal{N}^\sharp_{\mu, - k\mu}$ is equal to  the kernel of $\square U_{-1}  : \mathcal{N}_{\mu,  - k\mu} \to \mathcal{N}_{\mu, 1 - k\mu}$. 
\end{prop}

\parag{Proof} Recall that $K_\mu := D_N\big/H_\mu$ where $H_\mu$ is generated by $\partial_1, \dots, \partial_{k-1}, E+ \mu$. To show that $\chi$ exists it is enough to show that $\partial_h, h \in [1, k-1]$ and $E + \mu$ annihilate the class of  $U_{-1}$ in $\mathcal{N}_{\mu, 1-k\mu}$. This is given by formulas $E_h(\mu), h \in [1, k]$. \\
To prove the injectivity of $\chi$, simply remark that $K_\mu$ is simple and $U_{-1}$ is not in $\mathcal{J}_{\mu, 1- k\mu}$ thanks to Lemma \ref{2/3}.\\
The last assertion is easy since $\square (U_1 + \mu\sigma_1)$ vanishes on $K_\mu$,  the image of  $\square (U_1 + \mu\sigma_1)$ contains $\partial_1, \dots, \partial_{k-1}, E + \mu$,  the generators of $H_\mu$, thanks to the relations $F_j(\mu)$ for $j \in [1, k]$ (we assume that $k\mu \not= 1$)  and since the ideal $H_\mu$ is maximal (as $K_\mu$ is simple).$\hfill \blacksquare$\\

As a consequence of  Proposition \ref{magic 1} we have two exact sequences
\begin{equation*}
 0 \to \mathcal{N}^\sharp_{\mu, -k\mu} \to \mathcal{N}_{\mu, -k\mu} \to  K_\mu \to 0 \tag{$\lambda = -k\mu$},
 \end{equation*}
 and
\begin{equation*}
 0 \to K_\mu \to \mathcal{N}_{\mu, 1 - k\mu} \to \mathcal{N}^\square_{\mu, 1 - k\mu}   \to  0. \tag{$\lambda = 1 - k\mu$}
 \end{equation*}

where $K_\mu$ is a simple $D_N$-module since we assume that $\mu \not\in \mathbb{Z}$  and $k\mu \not= 1$.

\begin{thm} Assume that $\mu \not\in \mathbb{Z}$ and $k\mu \not= 1$. Then the $D_N$-modules $\mathcal{N}^\sharp_{\mu, -k\mu} $ and $\mathcal{N}^\square_{\mu, 1 - k\mu} $ are minimal extensions of trivial holomorphic rank $(k-1)$-vector bundles with  simple pole meromorphic connections along the hyper-surface $\{\Delta = 0 \}$.\\
Moreover, the map $\square(U_1 + \mu\sigma_1) : \mathcal{N}^\square_{\mu, 1 - k\mu} \to  \mathcal{N}^\sharp_{\mu, -k\mu}$ induces an isomorphism between them.
\end{thm}

\parag{Proof} The only point which is not already proved in Proposition \ref{magic 1} is the fact that the isomorphic $D_N$-modules are minimal extensions of trivial holomorphic rank $(k-1)$-vector bundles with  simple pole meromorphic connections along the hyper-surface $\{\Delta = 0 \}$. Since $E + \mu$ is in the sub-module generated by $\partial_1, \dots, \partial_{k-1}$ in $\mathcal{N}_{\mu, -k\mu}$ thanks to the equality
$$ k(E + \mu) = U_0 + k\mu + \sum_{h=1}^{k-1} (k-h)\sigma_h\partial_h, $$
we see that $\mathcal{N}^\sharp_{\mu, -k\mu}$ is generated by $\partial_1, \dots, \partial_{k-1}$ (using $k\mu \not= 0$). But $K_\mu$ is a rank $1$ trivial vector bundle with a simple pole meromorphic
connection and then $\mathcal{N}^\sharp_{\mu, -k\mu}$ is the trivial holomorphic vector bundle with basis $\partial_1, \dots, \partial_{k-1}$ inside $\mathcal{N}_{\mu, -k\mu}$.\\
Since $\mathcal{N}^\sharp_{\mu, -k\mu}$ has no torsion and $\mathcal{N}^\square_{\mu, 1-k\mu}$ has no co-torsion, these isomorphic $D_N$-modules are minimal extensions. $\hfill \blacksquare$\\

\begin{cor}\label{simple factors 2} The Theorem above  describes the structure of the $D_N$-modules $\mathcal{N}_{\mu, \lambda}$ when $\mu \not\in \mathbb{Z}$ assuming that $ k\mu \not= 1$   for $\lambda = n - k\mu$ with $n \in \mathbb{Z}$ thanks to Theorem \ref{isom 2}. $\hfill \blacksquare$\\
\end{cor}

\parag{Examples} For $k = 2$ the rank $1$ vector bundle associated to $\mathcal{N}^\sharp_{\mu, -2\mu}$ is trivial with basis $\partial_1$. It is isomorphic to $D_N\Delta^{-(\mu+1/2)}$.

\parag{Proof} We are computing the $D_N$-module $D_N\partial_1\big/ \mathcal{I}_\mu + D_N(U_0 + 2\mu)$ where the ideal $\mathcal{I}_\mu$ is equal to $ D_NT^2_\mu$.
 So we have the relations
\begin{align*}
& U_0 + 2\mu = \sigma_1\partial_1 + 2\sigma_2\partial_2 + 2\mu = 0 \\
& T^2_\mu =  \partial_1^2 + \partial_2(\sigma_1\partial_1 + \sigma_2\partial_2 + \mu) = 0 \\
& \partial_1(\sigma_1\partial_1 + 2\sigma_2\partial_2 + 2\mu ) = 0 \quad {\rm which \ implies} \\
& \sigma_1\partial_1^2 + 2\sigma_2\partial_1\partial_2 = -(2\mu + 1)\partial_1
\end{align*}
This gives, denoting by $e := \partial_1$ our generator,  since  $(E + 2\mu) \in  \mathcal{O}_N\partial_1$:
\begin{align*}
& 2\partial_1e + \sigma_1\partial_2e = 0 \\
& \sigma_1\partial_1e + 2\sigma_2\partial_2e = -(2\mu + 1)e
\end{align*}
The conclusion follows easily since $\Delta = \sigma_1^2 - 4\sigma_2$.$\hfill \blacksquare$\\

\subsection{The cases of  $\lambda = 0$ and $\lambda = -1$ with $\mu \not\in \mathbb{Z}$}

We define $\mathcal{N}^\sharp_{\mu, 0}$ as the kernel of the map $\mathcal{N}_{\mu, 0} \to \mathcal{O}_N$ sending $1$ to $1$. It is also the kernel of the map
$$ \square(U_1 + \mu\sigma_1) : \mathcal{N}_{\mu, 0}  \to \mathcal{N}_{\mu, -1} $$
thanks to formulas $F_j(\mu), j \in [1, k]$ and the fact that $U_1 + \mu\sigma_1$ is not in  $\mathcal{I}_\mu + D_N(U_0 + 1)$ (see lemma \ref{remplace} and the proof of lemma \ref{2/3}). \\
So  $\mathcal{N}^\sharp_{\mu, 0}$  is the quotient of the left ideal  in $D_N$  generated by  the ideal $\partial_1, \dots, \partial_k$ by the ideal  $I_\mu + D_NU_0$. So there is the relation
  between $\partial_1, \dots, \partial_k$ given by $U_0$ : $\sum_{h=1}^k h\sigma_h\partial_h = 0$.  On the open subset  $\sigma_k\Delta(\sigma) \not= 0$ we have a trivial rank $k-1$ holomorphic vector bundle with a meromorphic connection  having  a simple pole on  the hyper-surface $\{\sigma_k\Delta(\sigma) = 0\}$. Since $\mathcal{N}_{\mu, 0}$ has no torsion, $\mathcal{N}^\sharp_{\mu, 0}$ has no torsion.\\

We define $\mathcal{N}^\square_{\mu, -1}$ as the quotient of $\mathcal{N}_{\mu, -1}$ by  the kernel of the $D_N$-linear map  $\square U_{-1} : \mathcal{N}_{\mu, -1} \to \mathcal{N}_{\mu, 0}$ which contains 
 the image (isomorphic to $\mathcal{O}_N$)  of the map $\square(U_1 + \mu\sigma_1) : \mathcal{N}_{\mu, 0} \to \mathcal{N}_{\mu, -1}$ since this kernel contains  $(U_1 + \mu\sigma_1)$ by formula $(T)$.\\
So we have an induced  injective map $\square U_{-1} : \mathcal{N}^\square_{\mu, -1} \to \mathcal{N}_{\mu, 0} $ whose image in inside $\mathcal{N}^\sharp_{\mu, 0}$ since $U_{-1}$ belongs to the sub-module generated by $\partial_1, \dots, \partial_k$. But for $k\mu \not= 1$ the formulas $E_j(\mu)$ implies that the image of $\square U_{-1}$  contains $\partial_1, \dots, \partial_k$. So this map is surjective and 
 we conclude that  $\square U_{-1}$ is an isomorphism between  $\mathcal{N}^\square_{\mu, -1}$ and $\mathcal{N}^\sharp_{\mu, 0}$.\\

 Since $\mathcal{N}^\square_{\mu, -1}$ has no co-torsion and  $\mathcal{N}^\sharp_{\mu, 0}$  has no torsion, both $\mathcal{N}^\square_{\mu, -1}$ and $\mathcal{N}^\sharp_{\mu, 0}$ are minimal extensions of their respective restrictions to $\{\sigma_k\Delta(\sigma) \not= 0 \}$ which are rank $k-1$ coherent sheaves without torsion having  a simple pole meromorphic connection on  the hyper-surface $\{\sigma_k\Delta(\sigma) = 0\}$. \\
 The following lemma precise that, for $k\mu \not= 1$,  the inclusion of the image of the map $\square (U_1 + \mu\sigma_1)$ in $\mathcal{N}^\sharp_{\mu, -k\mu}$ is an isomorphism. As a corollary we obtain that $\mathcal{N}^\square_{\mu, -1}$ is a rank $(k-1)$ free  $ \mathcal{O}_N$-module with basis $\partial_1, \dots, \partial_{k-1}$ for $k\mu \not= 1$.
  
 \begin{lemma}\label{10/3} For $k\mu \not= 1$ the quotient  map from $ \mathcal{N}_{\mu, -1}\big/D_N(U_1 + \mu\sigma_1)$  to $ \mathcal{N}^\square_{\mu, -1}$ is an  isomorphism.
This implies the equality of the kernel of $\square U_{-1}$ with the image of $\square(U_1 + \mu\sigma_1)$ inside $\mathcal{N}_{\mu, -1}$. So $\mathcal{N}^\sharp_{\mu, -k\mu}$ is isomorphic to $\mathcal{O}_N$ as a $D_N$-module.
\end{lemma}

\parag{Proof} Let us compute the quotient $Q$  of $\mathcal{N}_{\mu, -1}\big/D_N(U_1 + \mu\sigma_1)$ by the sub-module generated by $\partial_1, \dots, \partial_{k-1}$. The elements of order at most 1 in $Q$ are generated on $\mathcal{O}_N$  by $1$ and $\partial_k$ with the relations
$$ k\sigma_k\partial_k + 1 = 0 \quad {\rm and} \quad \sigma_1\sigma_k\partial_k + \mu\sigma_1 = 0 .$$
Writing the second relation $\sigma_1(\sigma_k\partial_k + 1/k) + (\mu - 1/k)\sigma_1$ we find that $\sigma_1 = 0$ if $\mu \not= 1/k$. Then we obtain that $\sigma_1^2\partial_1 = 0$ and $\sigma_1\partial_h = 0 \quad \forall h \in [2, k]$ in $Q$ and then $Q$ is of $\sigma_1$-torsion. But this is not possible because the co-normal of  the hyper-surface $\{\sigma_1 = 0 \}$ is not in the characteristic variety of $Q$.
So $Q = 0$.$\hfill \blacksquare$\\

Since $D_N(U_1 + \mu\sigma_1) \subset \mathcal{N}_{\mu, -1}$ is isomorphic to $\mathcal{O}_N$ we obtain, for $k\mu \not= 1$, the exact sequence
$$0 \to \mathcal{O}_N \to  \mathcal{N}_{\mu, -1} \to \mathcal{N}_{\mu, -1}\big/D_N(U_1 + \mu\sigma_1) \to 0 $$
where the last term is associated to a rank $k-1$ trivial vector bundle (with $\mathcal{O}_N$-basis $\partial_1, \dots, \partial_{k-1}$)  with a simple pole meromorphic connection along the hyper-surface $\{\sigma_k\Delta(\sigma) = 0 \}$.

 So we obtain  the following result.
 
 \begin{thm}\label{8/5} For $k\mu \not= 1$ we have the exact sequences of $D_N$-modules
 \begin{align*}
 & 0 \to \mathcal{O}_N  \to   \mathcal{N}_{\mu, -1} \to \mathcal{N}^\square_{\mu, -1} \to 0 \\
 & 0 \to \mathcal{N}^\sharp_{\mu, 0} \to \mathcal{N}_{\mu, 0} \to \mathcal{O}_N \to 0
 \end{align*}
 where $\mathcal{O}_N$ is send to $ \mathcal{N}_{\mu, -1} $ by $1 \mapsto U_1 + \mu\sigma_1$. Moreover $\square U_{-1}$ induces an isomorphism of $\mathcal{N}^\square_{\mu, -1} $ to $\mathcal{N}^\sharp_{\mu, 0}$ and these two $D_N$-modules are minimal extensions of trivial rank $k-1$ holomorphic vector bundles  with a simple pole connection on  the hyper-surface $\{\sigma_k\Delta(\sigma) = 0\}$.$\hfill \blacksquare$\\
 \end{thm}

 \parag{An example} We detail the case $k = 2$ for  the $D_N$-module $\mathcal{N}^\sharp_{\mu, 0}$ which is the quotient  $(D_N\partial_1 + D_N\partial_2)\big/\mathcal{I}_\mu + D_NU_0$. It corresponds to the quotient   left ideal generated by
  
  \begin{align*}
& T^2_\mu = \partial^2_1 + \sigma_1\partial_1\partial_2 + \partial_2\big(\sigma_2\partial_2 + \mu\big) = 0 \\
& \sigma_1\partial_1 + 2\sigma_2\partial_2 = 0
\end{align*}
We obtain in $\mathcal{N}^\sharp_{\mu, 0}$
\begin{align*}
& \partial_1^2 + \frac{\sigma_1}{2}\partial_1\partial_2 + \mu\partial_2 = 0 \quad {\rm and \ after \ multiplication \ by} \ \sigma_2\\
& \sigma_2\partial_1^2 + \frac{\sigma_1}{2}\partial_1(-\frac{\sigma_1}{2})\partial_1 + \mu(-\frac{\sigma_1}{2})\partial_1 = 0   \quad {\rm and \ so}: \\
& \Delta\partial_1^2 = -(2\mu + 1)\sigma_1\partial_1
\end{align*}
Now we consider on $\{\sigma_2\Delta(\sigma) \not= 0\}$ the free rank $1$ $\mathcal{O}_N$-module with basis $e := \partial_1$ and we obtain
$$ \partial_1e =  -\frac{2\mu + 1}{\Delta}\sigma_1e .$$
The relation
$$ \partial_1\big(\sigma_1\partial_1 + 2\sigma_2\partial_2 \big) = 0$$
gives
\begin{align*}
& \partial_1 + \sigma_1\partial_1^2 + 2\sigma_2\partial_2\partial_1 = 0 \\
& 2\sigma_2\partial_2e = -e + \frac{2\mu+1}{\Delta}\sigma_1^2e \quad {\rm and \ using } \quad \sigma_1^2 = \Delta + 4\sigma_2 \\
& \partial_2e = \frac{\mu}{\sigma_2} + 4\frac{\mu + 1/2}{\Delta}e 
\end{align*}
Then the $D_N$-module $\mathcal{N}^\sharp_{\mu, 0}$ for $k=2$ is the minimal extension of the rank $1$ trivial holomorphic line bundle corresponding to function  $e := \sigma_2^\mu\Delta^{-(\mu + 1/2)}$.
Note that the weight of $e$ is $-1$ !\\
Note that for $\mu = 1/2$ we find $ e = \sigma_2^{1/2}\Delta^{-1}$

\parag{The case $k = 2, \mu = 1/2, \lambda = -1$ } We detail\footnote{For $\lambda = -1$ and $k = 2$ the only case where $\lambda = -k\mu$ is for $\mu = 1/2$.} for $k = 2$ and $\mu = 1/2$  the  quotient map form  the  $D_N$-module $H_\mu/(\mathcal{I}_\mu + D_N(U_0 + 1))$ to  $\mathcal{N}^\square_{\mu, -1}$ deduced from the inclusion of the image of $\square(U_1 + \mu\sigma_1)$ in the kernel of $\square U_{-1}$ in $\mathcal{N}_{1/2, -1}$. \\
 Recall that $H_\mu := Ker (\square U_{-1})$ is generated by $\partial_1, \dots, \partial_{k-1}$ and $E + \mu$ thanks to formulas $F_j(\mu), j \in [1, k]$, and that $(E + \mu)$  is not necessary for $\lambda = -k\mu$ since $k(E + \mu) - (U_0 -\lambda)$ is in $\partial_1, \dots, \partial_{k-1}$  in this case. Then $\mathcal{N}^\sharp_{1/2, -1}$ corresponds to the quotient of $H_\mu$ by the left ideal generated by
  
  \begin{align*}
& T^2_\mu = \partial^2_1 +  \partial_2\big(\sigma_1\partial_1 + \sigma_2\partial_2 + 1/2\big)  \quad {\rm and}  \\
&U_0 + 1 =  \sigma_1\partial_1 + 2\sigma_2\partial_2 + 1
\end{align*}
We obtain
\begin{align*}
& \partial_1^2 + \frac{\sigma_1}{2}\partial_1\partial_2  = 0 \quad {\rm and \ after \ multiplication \ by} \ \sigma_2\\
& \sigma_2\partial_1^2 + (-\frac{\sigma_1}{2}\partial_1 -1/2)\frac{\sigma_1}{2}\partial_1  \\
& \Delta\partial_1^2 = - 2\sigma_1\partial_1 
\end{align*}
Comparing
\begin{align*}
& \partial_1(\sigma_1\partial_1 + 2 \sigma_2\partial_2 + 1) = 0 \quad {\rm with} \\
& \partial_1^2 + \frac{\sigma_1}{2}\partial_1\partial_2 = 0 \quad {\rm gives}\\
& \Delta \partial_2\partial_1 = 4\partial_1
\end{align*}
Since the quotient is generated by $e := \partial_1$ since 
$$2(E + 1/2) = \sigma_1\partial_1 + (\sigma_1\partial_1 + 2\sigma_2\partial_2 + 1) = \sigma_1\partial_1$$
 in $\mathcal{N}^\sharp_{1/2, -1},$
we obtain that this $D_N$-module is isomorphic to the rank $1$ trivial vector bundle with basis $e$ and meromorphic connection defined by
$$ \partial_1e = -\frac{2\sigma_1}{\Delta}e \quad {\rm and} \quad  \partial_2 e = \frac{4}{\Delta}e $$
so $e = \Delta^{-1} $.\\
Note that, since $\lambda = -1$ the weight of $\partial_1 = e$ has to be $-2$.\\

\parag{The case $k = 2, \mu = 1/2, \lambda = -1, \lambda = 0$ }  We want to describe in this case the sub-module (isomorphic to $\mathcal{O}_N$) $D_N(U_1 + (1/2)\sigma_1)$ of the kernel of $\square U_{-1}$ in $\mathcal{N}_{1/2, -1}$  and the sub-module $K_{1/2}$ of the kernel of $\square (U_1 + (1/2)\sigma_1)$ in $\mathcal{N}_{1/2, 0}$ for $k = 2$.\\

This will show that for $k\mu = 1$ neither the equality between  the kernel  of $U_{-1}$ and the  image of  $\square(U_1 + \mu\sigma_1)$ in $\mathcal{N}_{1/k, -1}$ nor between the image of $\square U_{-1}$ and the kernel of $(U_1 + (1/k)\sigma_1)$ in $\mathcal{N}_{1/k, 0}$  are no longer true in general.\\

We have already see that for $k = 2$ the $D_N$-module  $\mathcal{N}^\sharp_{1/2, -1}$ is a free rank $1$ $\mathcal{O}_N$-module with basis $e = \partial_1$ and isomorphic to $D_N \sigma_2^{1/2}\Delta^{-1}$. 
The quotient of $\mathcal{N}^\sharp_{1/2, -1}$ by the image of $\square (U_1 + \sigma_1/2)$ is of $\Delta$-torsion since we have in $\mathcal{N}^\sharp_{1/2, -1}$
 $$ U_1 + \sigma_1/2 = \sigma_1(E + 1/2) - 2\sigma_2\partial_1 = \sigma_1^2\partial_1 + \sigma_1(-\frac{\sigma_1}{2}\partial_1 - 1/2) - 2\sigma_2\partial_1 = \frac{\Delta}{2}\partial_1.$$
 And the isomorphism $\mathcal{N}^\sharp_{1/2, -1} \to D_N\Delta^{-1}$ sending $e$ to $\Delta^{-1}$ send $(U_1 + \sigma_1/2)$ to $1/2$ realizing the isomorphism of the image of 
 $\square (U_1 + \sigma_1/2)$  to $\mathcal{O}_N$.\\
 
 Finally, the quotient of  $\mathcal{N}^\sharp_{1/2, 0}$ by the image of $\square U_{-1}$ which is isomorphic to $K_{1/2} \simeq D_N \sigma_2^{-1/2}$ is of $\sigma_2$-torsion since in  $\mathcal{N}^\sharp_{1/2, 0}$ we have
$$ U_{-1} = 2\partial_1 + \sigma_1\partial_2$$
 and so 
 $$ \sigma_2U_{-1} = 2\sigma_2\partial_1 + \sigma_1\sigma_2\partial_2 = 2\sigma_2\partial_1 + \sigma_1(\frac{-\sigma_1}{2}\partial_1) = -\frac{\Delta}{2}\partial_1 $$
 and the isomorphism $\mathcal{N}^\sharp_{1/2, 0} \to D_N\sigma_2^{1/2}\Delta^{-1}$ sends $U_{-1}$ to $D_N\sigma_2^{-1/2} \simeq K_{1/2}$.\\

\subsection{The map $\square(E + \mu)$ for $\lambda =  -k\mu -k+h, h \in [1, k]$}

To conclude, we examine the cases for which the map  $\square(E + \mu)$ is not an isomorphism.

\begin{prop}\label{29/3} For $\mu \not\in \mathbb{Z}$ and  $h = 1$ the class of $U_{-1}$\footnote{ which generates a sub-module isomorphic to $K_{\mu + 1}$.} in $\mathcal{N}_{\mu+1, -k(\mu+1) +1}$ is in the kernel of 
$\square(E + \mu) : \mathcal{N}_{\mu+1, -k(\mu +1)+1} \to \mathcal{N}_{\mu, -k\mu -k+1}$.\\
Moreover, for  $h = k$ the image of $\square(E + \mu): \mathcal{N}_{\mu+1, -k(\mu + 1) + k} \to \mathcal{N}_{\mu, -k\mu}$ contains $\partial_1, \dots, \partial_{k-1}$. So it is equal to 
$\mathcal{N}^\sharp_{\mu, -k\mu}$(see section 5.2).
\end{prop}

\parag{Proof} For $h = 1$ consider the formula $(E_1(\mu))$ and the commutation relation $(E + \mu)U_{-1} = U_{-1}(E + \mu) - k\partial_1$. This implies the relation
\begin{equation*}
-\partial_1(U_0 + k\mu + k - 1) + U_{-1}(E + \mu) \in \mathcal{I}_\mu 
\end{equation*}
which is enough to give the first assertion.\\
For $h = k$ the fact that this image contains $\partial_1, \dots, \partial_{k-1}$, thanks to the remark following Lemma \ref{4/9},  is enough to conclude since $\mathcal{N}^\sharp_{\mu, -k\mu}$ is maximal because $K_\mu$ is simple for $\mu \not\in \mathbb{Z}$.$\hfill \blacksquare$\\

\begin{lemma}\label{28/4} We have the commutation relation
\begin{equation}
(U_1 + (\mu +1)\sigma_1)(E + \mu) = (E+\mu)(U_1 + \mu\sigma_1) 
\end{equation}
\end{lemma}

\parag{Proof} Recall that $U_1 + \mu\sigma_1 = \sigma_1(E + \mu) - \sum_{h=1}^{k-1} (h+1)\sigma_{h+1}\partial_h $. We have
\begin{align*}
& [E, \sigma_1] = \sigma_1 \quad  [E, \sigma_{h+1}\partial_h] =  0 \   \forall h \in [1, k-1] \quad {\rm and \ so} \\
&   [E, U_1] = [E, \sigma_1E] = [E, \sigma_1]E = \sigma_1E. \ {\rm Then \ we \ obtain} \\
& [U_1 + \mu\sigma_1, E + \mu] = -\sigma_1E + \mu[\sigma_1, E] = -\sigma_1(E+ \mu).
\end{align*}
and the formulas $(30)$ follows.$\hfill \blacksquare$\\

\begin{cor}\label{28/4 bis}
So we obtain  for each $\mu$ and $\lambda$ the commutative diagram
$$ \xymatrix{ \mathcal{N}_{\mu+1, \lambda}  \qquad\ar[r]^{\square(U_1 + (\mu+1)\sigma_1)} \ar[d]_{\square(E + \mu)}   & \qquad  \mathcal{N}_{\mu+1, \lambda - 1} \ar[d]^{\square(E + \mu)} \\
\mathcal{N}_{\mu, \lambda} \qquad \ar[r]^{\square(U_1 + \mu\sigma_1)}  & \qquad \mathcal{N}_{\mu, \lambda - 1}}  $$
\end{cor}

\parag{Conclusion}  Proposition  \ref{29/3} and Corollary \ref{28/4 bis} shows that for  $\mu \not\in \mathbb{Z}$ and for $h \in [1, k]$ the kernel of $\square(E + \mu) : \mathcal{N}_{\mu + 1, -k\mu -k + h} \to \mathcal{N}_{\mu, -k\mu -k + h} $
is deduced from the isomorphisms $\square(U_1 + (\mu+1)\sigma_1) : \mathcal{N}_{\mu + 1, -k\mu -k + h} \to \mathcal{N}_{\mu + 1, -k\mu -k + h-1}$ and its image from the isomorphisms
$\square(U_1 + \mu\sigma_1) : \mathcal{N}_{\mu, -k\mu -k + h} \to \mathcal{N}_{\mu, -k\mu -k + h-1}$. But we know that for $h = 1$ the kernel  is generated by $U_{-1}$ so isomorphic to $K_{\mu + 1}$  and that for $h = k$ the image is equal to $\mathcal{N}^\sharp_{\mu, -k\mu}$.\\
This determines the kernel and the image of $\square(E + \mu): \mathcal{N}_{\mu+1, \lambda} \to \mathcal{N}_{\mu, \lambda}$ for each value of  $\lambda = -k\mu - h, h \in [0, k-1],$ up to these isomorphisms.\\

\end{document}